\documentclass{article}

\usepackage[T1]{fontenc}
\usepackage[utf8]{inputenc}
\usepackage[margin=1in,top=1in,bottom=3cm]{geometry}
\usepackage{perpage}
\MakePerPage{footnote}
\usepackage{epsfig}
\usepackage{latexsym}
\usepackage{amsmath}
\usepackage{mathrsfs}
\usepackage{yfonts}
\usepackage{makeidx}
\usepackage{graphicx}
\usepackage{slashed}
\usepackage{hyperref}
\usepackage[titletoc]{appendix}
\usepackage{chngcntr}
\usepackage{bm}
\usepackage{multicol}
\usepackage{nameref}
\usepackage{textcomp}
\usepackage{changepage}
\usepackage{lettrine}
\usepackage{enumitem}
\usepackage[sortcites,backend=bibtex,style=numeric,sorting=anyvt,doi=false,url=false, isbn=false]{biblatex}
\AtEveryBibitem{\clearfield{month}}
\AtEveryBibitem{\clearfield{day}}
\hypersetup{citecolor=blue,colorlinks=true,linkcolor=blue,pdfstartview={FitH},linktoc=page,pdftitle={Maxwell-Decay-on-DSRN-BH.pdf}pdfauthor={Mokdad Mokdad}}
\usepackage{amsthm,amssymb}
\usepackage{subcaption}
\newtheorem{thm1}{Theorem}%[section]
\newtheorem{prop1}[thm1]{Proposition}
\newtheorem{lem1}[thm1]{Lemma}

 \newtheoremstyle{TheoremNum}
        {\topsep}{\topsep}              %%% space between body and thm
        {\itshape}                      %%% Thm body font
        {}                              %%% Indent amount (empty = no indent)
        {\bfseries}                     %%% Thm head font
        {.}                             %%% Punctuation after thm head
        { }                             %%% Space after thm head
        {\thmname{#1}\thmnote{ \bfseries #3}}%%% Thm head spec
    \theoremstyle{TheoremNum}

     \newtheoremstyle{TheoremNum}
        {\topsep}{\topsep}              %%% space between body and thm
        {\itshape}                      %%% Thm body font
        {}                              %%% Indent amount (empty = no indent)
        {\bfseries}                     %%% Thm head font
        {.}                             %%% Punctuation after thm head
        { }                             %%% Space after thm head
        {\thmname{#1}\thmnote{ \bfseries #3}\thmnumber{}}%%% Thm head spec
    \theoremstyle{TheoremNum}

\renewcommand{\d}{\mathrm{d}}
\newcommand{\scri}{{\mathscr I}}

\newcommand{\R}{\mathbb{R}}
\newcommand{\dl}{\partial}

\newcommand{\hf}{\frac{1}{2}}

\makeindex

\title{Reissner-Nordstrøm-de Sitter Manifold : Photon Sphere and Maximal Analytic Extension}
\author{Mokdad Mokdad\thanks{email: mokdad.al.mokdad@gmail.com - Mokdad.Mokdad@univ-brest.fr} 
		\\LMBA -- Université de Bretagne Occidentale}

\begin{document}

\maketitle

\begin{abstract}
	{%We study the Reissner-Nordstr{\o}m-de Sitter spacetime which is a spherically symmetric solution of the Einstein-Maxwell system. It models a charged, non-rotating, spherically symmetric black hole spacetime in the presence of a positive cosmological constant. 
	This paper is devoted to the study of the Reissner-Nordstr{\o}m-de Sitter black holes and their maximal analytic extensions. We study some of their properties that lays the groundwork for obtaining (in separate papers) decay results \cite{mokdad_decay_2016} and constructing conformal scattering theories for test fields on such spacetimes \cite{mokdad_conformal_2016}. Here, we find the necessary and sufficient conditions on the parameters of the Reissner-Nordstr{\o}m-de Sitter metric ---namely, the mass , the charge, and the cosmological constant--- to have three horizons. Under this conditions, we prove that there is only one photon sphere and we locate it. We then give a detailed construction of the maximal analytic extension of the Reissner-Nordstr{\o}m-de Sitter manifold in the case of three horizons. 
		%%Up to our knowledge, this has never been carried out in the literature.
	}
\end{abstract}

\tableofcontents
 \newpage
\section{Introduction}

The year 2015 marked the 100th anniversary of Albert Einstein's presentation of the complete theory of General Relativity to the Prussian Academy. A hundred years have passed and Einstein's general theory of relativity is still the most accurate description of gravity that we ever had. According to this theory, gravity is the manifestation of the curvature of spacetime, which is a Lorentzian 4-manifold consisting of all the events in ``space'' and ``time'', where these two concepts merge into one. The field equations that govern the laws of gravity relate the presence of energy and momentum to the curvature of a Lorentzian metric which is a solution of the equations. The  tensorial form of the equations is
$$\mathbf{G_{ab}} + \Lambda\mathbf{ g_{ab}}= \frac{8 \pi G}{c^4} \mathbf{T_{ab}} \; .$$ 
The unknown in the equations is the Lorentzian metric $\mathbf{ g_{ab}}$ which is a non-degenerate symmetric (0,2)-tensor of signature $(+,-,-,-)$\footnote{Or $(-,+,+,+)$, the difference is a matter of taste in most situations. In this work, we shall carry on with the convention in the text above.}. The Einstein tensor $\mathbf{G_{ab}}$ is
$$\mathbf{G_{ab}}=\mathbf{R_{ab}}-\hf\mathbf{g_{ab}R}\; ,$$
where $\mathbf{R_{ab}}$ is the Ricci curvature tensor of the metric $\mathbf{ g_{ab}}$ and $\mathbf{R}$ is the scalar curvature of the metric. These curvature quantities are given by the Riemann curvature tensor $\mathbf{R_{abcd}}$ which itself is locally given in terms of the Christoffel symbols of the metric:
$$\mathbf{\Gamma^c{}_{ab}}=\hf \mathbf{g^{cd}}\left(\dl_\mathbf{a} \mathbf{g_{bd}}+\dl_\mathbf{b} \mathbf{g_{ad}}-\dl_\mathbf{d} \mathbf{g_{ab}} \right)\; ,$$
and
$$\mathbf{R^a{}_{bcd}} = \partial_{\mathbf{c}}\mathbf{\Gamma^a{}_{db}}
- \partial_{\mathbf{d}}\mathbf{\Gamma^a{}_{cb}}
+ \mathbf{\Gamma^a{}_{ce}\Gamma^e{}_{db}}
- \mathbf{\Gamma^a{}_{de}\Gamma^e{}_{cb}}\; .$$
The scalar curvature is the trace of the Ricci curvature which in turn is given by the trace of the Riemann curvature tensor: 
$$\mathbf{R_{ab}}=\mathbf{R^c{}_{acb}} \quad \textrm{ and } \quad \mathbf{R}=\mathbf{R^a{}_a}\; .$$    

$\mathbf{T_{ab}}$ is the {energy-momentum tensor}\footnote{Also called stress-energy-momentum tensor or {stress-energy tensor}\index{Stress-energy tensor}.} determined by the matter, energy, and momentum, present in the spacetime. The rest are constants: $\Lambda$ is the cosmological constant, $G$ is the gravitational constant of Newton, and finally $c$ is the speed of light in vacuum.  

\subsection*{{Black Holes}}

\begin{figure}
	\centering
	\includegraphics[scale=1]{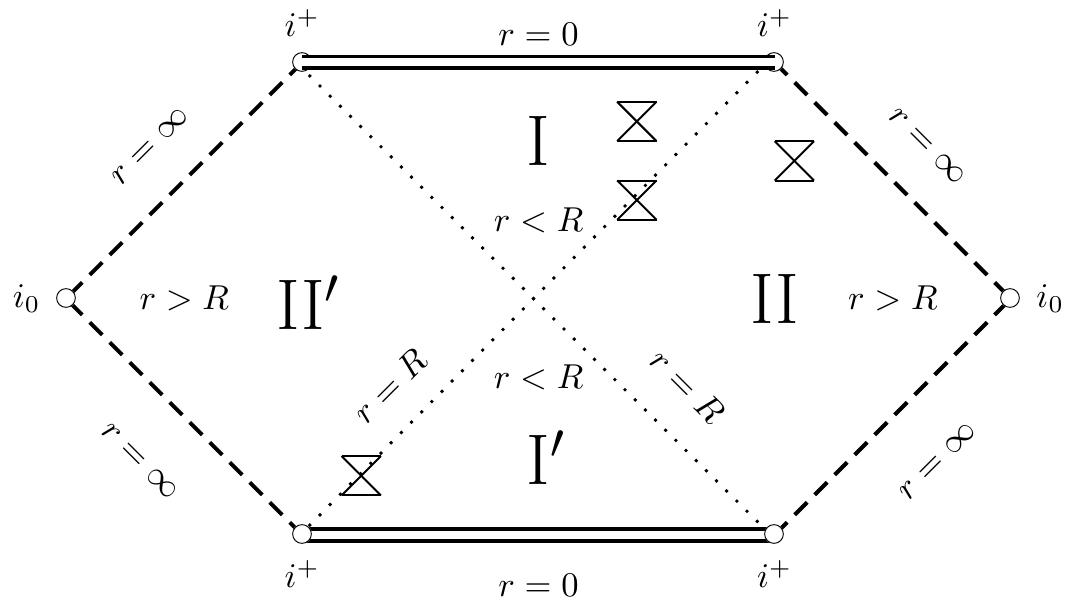}
	\caption{\emph{Penrose-Carter conformal diagram of Kruskal-Szekeres spacetime, the maximal extension of a Schwarzschild black hole (I and II).}}
	\label{fig:introSchwrzchildmax}
\end{figure}

The curvature terms in Einstein's field equations contain first and second order partial derivatives of the metric and they are a highly nonlinear system of partial differential equations, which makes them very hard to solve in general. However, several families of exact solutions are known. The trivial solution in the vacuum case, i.e. when the energy-momentum tensor vanishes, with a zero cosmological constant is the simplest Lorentzian metric on $\R^4$, 
$$\mathbf{g}=\d t^2 - \d x^2 - \d y^2 - \d z^2\; ,$$ which is known as the Minkowski metric. Minkowski spacetime is flat, meaning that the Riemann curvature tensor vanishes identically. The second best-known solution is the Schwarzschild metric. This solution of the Einstein's vacuum equations with zero cosmological constant, describes an empty spacetime outside a non-rotating and uncharged spherical body of mass $M$ and radius $R=2MGc^{-2}$ by a metric $\mathbf{g}$ defined on $\mathcal{M}=\R_t \times ]0,+\infty[_r \times \mathcal{S}^2_{\theta,\varphi}$, and whose spherical coordinate expression is 
$$\mathbf{g}=\left(1-\frac{R}{r}\right)\d t^2 - \left(1-\frac{R}{r}\right)^{-1}\d r^2 -r^2 \left(\d \theta^2 + \sin(\theta)^2\d \varphi^2\right)\; .$$ 
The metric describes the spacetime region with $r>R$, called the exterior Schwarzschild solution. Nevertheless, one could (mathematically at least) assume that the space inside the region where the body is supposed to be, is another empty region of spacetime given by the same metric expression but for $0<r<R$, which is called the interior Schwarzschild solution. When viewed with these coordinates, the two solutions appear as two completely separate solutions with no physical connection between them, separated by an apparent singularity at $r=R$. However, viewed as a Lorentzian manifold, the singularity at $r=R$ is a mere coordinate singularity due to this particular choice of coordinates. In fact, the Kruskal-Szekeres coordinates extends the original Schwarzschild spacetime and cover the entirety of $r>0$. The Kruskal-Szekeres extension is the maximal analytic extension of the Schwarzschild spacetime, and it describes a theoretical eternal black/white hole. It shows that the hypersurface at $r=R$ is not singular but rather a regular null hypersurface that acts like a barrier which can be crossed only in one direction and is therefore an event horizon, hence the name black/white hole. (Figure \ref{fig:introSchwrzchildmax})

The singularity at $r=0$ is different. This is a genuine physical or geometrical singularity since the scalar curvature
$$\mathbf{R}=\frac{12R^2}{r^6}$$
clearly blows up, and since this is a scalar quantity, it means no coordinate transformation could resolve the singularity at $r=0$.

One form of energy that can induce gravity is light, i.e. electromagnetic radiation. A source-free Maxwell field on spacetime is a 2-form $F$ satisfying Maxwell's equations:
$$\d F=0 \qquad;\qquad \d \star F=0$$ 
where $\d$ is the exterior differentiation and $\star$ is the Hodge star operator. The Maxwell system describes the phenomena of electromagnetism. The presence of a Maxwell field curves spacetime around it. Its effect is given by the electromagnetic energy-momentum tensor
$$\mathbf{T_{ab}}=\frac{1}{4} \mathbf{g_{ab}}F^{cd}F_{cd}-F_{ac}{F_b}^c \; .$$

The Einstein-Maxwell equations are Einstein's field equations with $\mathbf{T_{ab}}$ given by the electromagnetic energy-momentum tensor of a Maxwell field. A known family of exact solutions to the Einstein-Maxwell equations with no cosmological constant generalizes the Schwarzschild solution by describing the spacetime outside a rotating charged black hole.  This family of solutions is known as the Kerr-Newman black holes, and each member of the family is uniquely described by three real parameters of the black hole: a mass $M$, a charge $Q$, and an angular momentum $a$. The Schwarzschild solution corresponds to the case where $a$ and $Q$ are both zero. If only $a=0$, the solution is called Reissner-Nordstr{\o}m black hole and it describes a non-rotating but charged black hole. When $Q=0$ but $a\ne 0$, the solution represents a rotating black hole, also known as a Kerr black hole. When a cosmological constant is present, a solution of the Einstein-Maxwell equations will have what can be called a de Sitter ($\Lambda>0$) or an anti-de Sitter ($\Lambda<0$) ``aspect''. The simplest solution to Einstein's vacuum equations with a (negative) positive  cosmological constant is the (anti-) de Sitter spacetime. The de Sitter spacetime is the analogue in Minkowski spacetime, of a sphere in ordinary Euclidean space. It is maximally symmetric, has constant positive scalar curvature, and is simply connected. It can be visualized as hyperboloid in a 5-dimensional flat Lorentzian manifold. Spherically symmetric asymptotically de Sitter spacetimes that are solutions to  the Einstein-Maxwell equations have metrics similar to the Schwarzschild metric. For example, in spherical coordinates, the metric of such a spacetime is typically  of the form:
$$\mathbf{g}=f(r)\d t^2 - \frac{1}{f(r)}\d r^2 -r^2 \left(\d \theta^2 + \sin(\theta)^2\d \varphi^2\right)\; ,$$
with $f$ possibly of the form
\begin{equation*}
f(r)=1-\frac{2MG}{c^2r}+\frac{Q^2}{r^2}-\Lambda r^2.\tag{{\small{I}}}\label{finintrothesis}
\end{equation*}

With all of $M,Q,$ and $\Lambda$ equal to zero, i.e. $f(r)=1$, we get Minkowski spacetime. When all but $\Lambda$ are equal to zero, we get the (anti-) de Sitter spacetime. When all but $M$ are zero, this is Schwarzschild's black hole. And when only $Q$ is zero we have  the Schwarzschild-de Sitter spacetime. Alternatively, if only $\Lambda$ equals zero, the spacetime is a Reissner-Nordstr{\o}m black hole. Finally, if all three parameters are non zero, we get the \emph{Reissner-Nordstr{\o}m-de Sitter black hole}, which in this work, is the spacetime we are interested in.

{In the rest of the paper, the metric is presented in units where both $G$ and $c$ are 1. Furthermore, we assume\footnote{The positivity of the cosmological constant $\Lambda$ is motivated by its experimental value. As far as we know, it has a very small value (~$10^{-122}$ \cite{barrow_value_2011}) but yet positive.} $M,\Lambda>0$  and $Q\ne 0$, and there is no rotation $(a=0)$}. So to fix notations, the Reissner-Nordstr{\o}m-de Sitter metric we are studying is given in spherical coordinates by 
\begin{equation}\label{RNdSmetric}
g_\mathcal{M}=f(r)\d t^2-\frac{1}{f(r)}\texttt{d}r^2-r^2\d \omega^2,
\end{equation}
where
\begin{equation}\label{f(r)}
f(r)=1-\frac{2M}{r}+\frac{Q^2}{r^2}-\Lambda r^2 \; ,
\end{equation}
and $\d \omega^2$ is the Euclidean metric on the $2$-Sphere, $\mathcal{S}^2$, which in spherical coordinates is
\begin{equation*}
\d \omega^2=\d \theta^2 + \sin(\theta)^2\d \varphi^2 \; ,
\end{equation*}
and $g_\mathcal{M}$ is defined on \index{$\mathcal{M}$}$\mathcal{M}=\R_t \times ]0,+\infty[_r \times \mathcal{S}^2_{\theta,\varphi}$ .
The Reissner-Nordstr{\o}m-de Sitter solution (which we sometimes abbreviate as ``\textbf{RNdS}''), is one of the spherically symmetric solutions of the Einstein-Maxwell field equations in the presence of a positive cosmological constant $\Lambda$. This solution models a non-rotating spherically symmetric charged black hole with mass $M$ and a charge $Q$, in a de Sitter background. The de Sitter background means that there is a cosmological horizon beyond which lies a region that stretches to infinity. While the Reissner-Nordstr{\o}m nature entails that near the singularity, depending on the relation between the mass and the charge, one has a succession of static\footnote{A spacetime, or part of it, is said to be static if the metric admits a timelike Killing vector field that is orthogonal to a family of spacelike hypersurfaces. It is said to be dynamic if there is no timelike Killing vector field.} and dynamic regions separated by the apparently singular hypersurfaces. The metric in these coordinates appears to have singularities at $r=0$ and at the zeros of $f$. Only the singularity at $r=0$ is a real geometric singularity at which the curvature blows up. The apparent singularities at the zeros of $f$ are artificial and due to this particular choice of coordinates. The regions of spacetime where $f$ vanishes are essential features of the geometry of the black hole, they are called event horizons or {horizons}\index{Horizon} for short, and $f$ is called the {horizon function}\index{Horizon function}. If $f$ has three positive zeros and a negative one, then the zeros in the positive range corresponds in an increasing order respectively to the {Cauchy horizon} or {inner horizon}, the {horizon of the black hole} or the {outer horizon}, and the {cosmological horizon. In our case, we work with three horizons corresponding to $r$ equals to $r_1,r_2$, and $r_3$, the three positive zeros of $f$. In this case, the region corresponding to $]0,r_1[$ is a static region in the interior of the black hole, and the one corresponding to $]r_2,r_3[$ is another static region in the exterior. %which we call $\mathcal{N}=\R_t \times ]r_2,r_3[_r\times\mathcal{S}^2_\omega$, 
An interior dynamic region separating the two static regions lies in $]r_1,r_2[$, and the region given by $r>r_3$ is a dynamic region near infinity. 
%Here $M$ is the mass of the black hole, $Q$ is its charge, and $\Lambda$ is the cosmological constant. We assume that $Q$ is real and non zero, and $M$ and $\Lambda$ are positive.

One can refer to classical books such as \cite{hawking_large_1973, wald_general_2010,chandrasekhar_mathematical_1984} for more on exact solutions and on Einstein's general theory of relativity. However, few works can be found on the RNdS spacetime. Some of the early works \cite{lake_reissner-nordstrom-sitter_1979, laue_maximally_1977,bronnikov_inverted_1979} shortly discuss the construction of maximally extended RNDS spacetimes. Works on global spacetime solutions and on their constructions that include RNdS cases can also be found for example in \cite{bolokhov_magnetic_2012,katanaev_geometric_1996,katanaev_global_1999}. Thus, the general aspects of the construction of the maximal analytic extension of RNdS spacetime is in the literature, but up to our knowledge, it has never been explicitly carried out in details. In part of this paper, we try to fill the ``gap'' by giving conditions on the free parameters of the RNDS metric and by constructing the maximal analytic extension of this spacetime in the most complete case (three horizons) in rather sufficient details, discussing the role of radial null geodesics in specifying the coordinates used and in obtaining the different extensions. We also discuss some of the geometrical and causal properties of these extensions (in the spirit of \cite{hawking_large_1973}).

On the other hand, some of the geometrical aspects of this spacetime are important for us in other works \cite{mokdad_decay_2016, mokdad_conformal_2016}. In these works, we are interested in the decay in time of Maxwell fields in the static region between the horizon of the black hole and the cosmological horizon, i.e. the exterior static region, of the RNdS black hole spacetime. The decay of test fields (such as the electromagnetic fields) plays an important role in studying the stability of solutions of Einstein's equation. The RNdS black holes can be considered as spherically symmetric models of the more important Kerr family of black holes which are believed to best represent real black holes that may be existing now in our universe. The exterior static region of the RNdS spacetime contains a photon sphere, i.e. null geodesics orbiting the black hole at fixed $r$. It is known that this phenomenon is an important part of black hole spacetimes geometry. Discussions of photon orbits and there effects can be found in \cite{khoo_lux_2016,chakraborty_aspects_2015,chakraborty_equilibrium_2015,chakraborty_strong_2016} for example. The effect of photon sphere concerning decay of test fields can be seen in  \cite{blue_decay_2008,andersson_hidden_2015,mokdad_decay_2016} among others. A priori, the existence of a photon sphere is an obstacle for the decay. Still, fields can decay, as shown in \cite{mokdad_decay_2016} or \cite{blue_decay_2008} and other works, but the photon sphere slows the decay as there will be null geodesics that rotate around the black hole near the photon sphere for arbitrary amount of time. We are also interested in constructing a conformal scattering theory on the exterior static region \cite{mokdad_conformal_2016}. For this, we need to have access to the boundary of the region corresponding to infinite $t$-values in the static interior region. This boundary is part of the maximal analytic extension of the spacetime which we construct in this paper.

This paper has two main sections:

\paragraph{Section \ref{sectionsetupandphotonsphere}:} We start the section by presenting the necessary and sufficient conditions (\ref{GC}) on the parameters $M,Q,$ and $\Lambda$ of the RNdS metric so that it has three horizons. We then verify our claim regarding these conditions along with the fact that there is a photon sphere only at one value of $r>0$ and it is located in the exterior static region. This is Proposition \ref{1photonsphere}, and up to our knowledge, this is not in the literature.

\paragraph{Section \ref{sec:maximalextension}:} This section is a detailed discussion and construction of the maximal analytic extension of the RNdS manifold in the case of three horizons at $0<r_1<r_2<r_3$. We start by exploring some properties of the black hole in the RNdS coordinates $(t,r,\omega)$. We then discuss the Regge-Wheeler $r_*$ coordinate and use it to obtain coordinate expressions of the radial null geodesics. Using the radial null geodesics we define the Eddington-Finkelstein advanced and retarded coordinates and extensions, showing that the event horizons are not singular but in fact are regular null hypersurfaces for the extended metric. The place where horizons of the same $r$-value ``meet'' is asymptotic to all of the Eddington-Finkelstein charts, these are the bifurcation spheres. To cover these spheres we need the Kruskal-Szekeres extensions. Each of these new extensions now cover all the horizons at $r=r_i$ and the bifurcation sphere where they intersect. Finally, we use the Kruskal-Szekeres charts to cover the manifold  of the maximal analytic extension. We discuss its causal structure, and some properties of its timelike singularity at $r=0$.

\section{Photon Sphere}\label{sectionsetupandphotonsphere}

In this section we study the horizon function $f$ given in (\ref{f(r)}) of the RNdS metric (\ref{RNdSmetric}). We put,
\begin{gather}\label{naming}
R=\frac{1}{\sqrt{6\Lambda}} \quad;\quad \Delta=1-12Q^2\Lambda \quad ; \quad m_1=R\sqrt{1-\sqrt{\Delta}} \quad ;\quad m_2=R\sqrt{1+\sqrt{\Delta}}\\ M_1=m_1-2\Lambda m_1^3 \quad ; \quad M_2=m_2-2\Lambda m_2^3 \label{naming2}\; . 
\end{gather}
and we consider the following conditions,
\begin{equation}\label{GC}
Q\neq0 \quad \textrm{and} \quad 0<\Lambda < \frac{1}{12Q^2} \quad \textrm{and} \quad M_1 < M < M_2 \; .
\end{equation}
The main result of this section is:
\begin{prop1}[Three Positive Zeros and One Photon Sphere]\label{1photonsphere}
	The function $f$ has exactly three positive distinct zeros if and only if (\ref{GC}) holds. In this case, there is exactly one photon sphere in the static exterior region of the black hole defined by the portion between the largest two zeros of $f$.
\end{prop1}

The proof is divided into parts: First, we study the conditions on $M,Q,$ and $\Lambda$ for $f$ to have three positive zeros, and then we show that in that case there is only one photon sphere.

\subsection{The Zeros of the Horizon Function}\label{Sec:Zerosoff}

The zeros of the function $f$ are the roots of the polynomial
\begin{equation}\label{P(r)}
r^2f(r)=P(r)=-\Lambda r^4 + r^2 -2Mr + Q^2 \, .
\end{equation}
Let us show that $P$ has exactly three positive and one negative real roots if and only if (\ref{GC}) holds. We will proof this in two lemmata.
\begin{lem1}\label{lem1}
	The polynomial $P$ has three positive roots if and only if
	\begin{equation}\label{Variation}
	P'(R)>0 \quad \textrm{and} \quad P(s_1)<0 \quad \textrm{and} \quad P(s_2)>0 \; ,
	\end{equation}
	
	where $0<s_1<s_2$ are the two positive roots of $P'$.
\end{lem1}

\begin{proof}
	The expressions of $P'$ and $P''$ are
	
	$$P'(r)=-4\Lambda r^3 + 2r -2M \;, \; P''(r)=-12\Lambda r^2 +2 \; ,$$
	
	and so $P''(R)=0$ . Because $R$ is the only positive root of $P''(r)$ and $P''(0)=2$, $P'$ is increasing on $[0,R]$ and decreasing on $[R,+\infty[$ with a local maximum at $R$. If $P'(R)$ is non positive, and since $P'(0)=-2M<0$, then $P'$ is everywhere non positive on $[0,+\infty[$ . Thus, $P$ is decreasing on $[0,+\infty[$ , and has only one root there as it decreases from $P(0)=Q^2>0$ to $-\infty$. Therefore, a necessary condition for $P$ to have three positive roots is that $P'(R)$ be positive. Clearly,
	\begin{equation}\label{P'(R)>0}
	P'(R)>0 \quad \Leftrightarrow \quad M<\frac{2}{3}R \; .
	\end{equation}
	As $P'(0)<0$, and $\lim_{r\rightarrow \pm\infty}P'(r)=\mp\infty$, then having a positive local maximum at $R$ implies that $P'$ has exactly two roots $0<s_1<R<s_2$ on the positive axis, and one on the negative axis. Also $P'$ changes sign after passing through each of its roots $s_1$ and $s_2$, which means that $P(s_1)$ and $P(s_2)$ are respectively the local minimum and the local maximum of $P$ over the interval $[0,+\infty[$ . We can conclude the following:
	
	\begin{itemize}
		\item If $P(s_1)>0$, then $P$ has one positive root $x$, with $s_2<x$.
		\item If $P(s_1)=0$, then $P$ has two positive roots $s_1$ and $x$, with $s_1<s_2<x$.
		\item If $P(s_1)<0$, then :
		\begin{itemize}
			\item If $P(s_2)<0$, then $P$ has one positive root $x$, with $x<s_1$.
			\item If $P(s_2)=0$, then $P$ has two positive roots $x$ and $s_2$, with $0<x<s_1<s_2$.
			\item If $P(s_2)>0$, then $P$ has three positive roots $r_1, r_2 ,$ and $r_3$, with $0<r_1<s_1<r_2<s_2<r_3$.
		\end{itemize}
	\end{itemize}
	This concludes the proof.
\end{proof}

Instead of finding $s_1$ and $s_2$ explicitly, we will, using the next lemma, transform the conditions in (\ref{Variation}) to those in (\ref{GC}) directly.

\begin{lem1}\label{lem2}
	If $P'(R)>0$, with $s_1$ and $s_2$ the positive roots of $P'$, then
	\begin{eqnarray}
	P(s_1)<0 &\textrm{if and only if}& P'(m_1)<0 \; ; \label{Betterconditions1}\\
	P(s_2)>0 &\textrm{if and only if}& P'(m_2)>0 \; , \label{Betterconditions2}
	\end{eqnarray}
	where $m_1$ and $m_2$ are defined in (\ref{naming}).
\end{lem1}

\begin{proof}
	We first note that
	\begin{eqnarray*}
		P(r) &=& -\Lambda r^4 + r^2 -2Mr + Q^2 \\
		&=& rP'(r)+T(r)
	\end{eqnarray*}
	where T is the polynomial
	
	$$T(r)=3\Lambda r^4 -r^2 + Q^2 \; .$$
	
	So, $P(s_1)=T(s_1)$ and $P(s_2)=T(s_2)$. Therefore if we study the sign of $T$ we shall know the sign of $P(s_1)$ and $P(s_2)$. Let $\bar{T}(r^2)=T(r)$, i.e. $$\bar{T}(r)=3\Lambda r^2 -r + Q^2 \; ,$$
	which has discriminant $\Delta=1-12\Lambda Q^2$. We investigate the different cases.
	\begin{itemize}
		\item If $\Delta<0$ then $\bar{T}$ has no real roots and is always positive, and hence so is $T$. In particular, this means that $T(s_1)$ and $T(s_2)$ are both positive, which is not the desired case.
		\item If $\Delta =0$ then $R$ is a double root for $\bar{T}$ and it is non negative. It follows that $T$ is also non negative, and the conditions of (\ref{Variation}) cannot be satisfied.
		\item Finally, if $\Delta >0$ which is $\Lambda<\frac{1}{12Q^2}$, then $\bar{T}$ has two positive roots $m_1^2 \; , \; m_2^2$ and hence $\pm m_1 \; , \; \pm m_2$ are the roots of $T$, and $T$ is positive on $[0,m_1[$, negative on $]m_1,m_2[$, and positive on $]m_2,+\infty[$ .
	\end{itemize}
	
	Thus, noting that $s_1$ and $m_1$ are strictly less than $R$, and $s_2$ and $m_2$ are strictly greater than $R$ when $\Lambda<\frac{1}{12Q^2}$ and assuming (\ref{P'(R)>0}), we see that
	\begin{eqnarray}\label{betterconditions1}
	P(s_1)=T(s_1)<0 &\textrm{if and only if}& m_1<s_1 \\
	P(s_2)=T(s_2)>0 &\textrm{if and only if}& m_2<s_2 \; . \label{betterconditions2}
	\end{eqnarray}
	
	Now the key point which makes the right hand sides of (\ref{betterconditions1}) and (\ref{betterconditions2}) more useful is that $P'$ is strictly monotonic on each side of $R$, and that $P'(s_1)=P'(s_2)=0$ . By applying $P'$ to (\ref{betterconditions1}) and (\ref{betterconditions2}), one gets (\ref{Betterconditions1}) and (\ref{Betterconditions2}).
	
\end{proof}

Let us summarize: Recalling $M_1$ and $M_2$ from (\ref{naming2}) and noting that $P'(m_1)=-4\Lambda R^3(1-\sqrt{\Delta})$, we see that $P$ has three positive roots if and only if
\begin{enumerate}
	\item $Q\neq 0$ and
	\item $0<\Lambda<\frac{1}{12Q^2}$ and
	\item $P'(R)>0 \quad i.e. \quad M<\frac{2}{3}R$ and
	\item $P'(m_1)<0 \quad i.e. \quad M_1<M$ and
	\item $P'(m_2)>0 \quad i.e. \quad M<M_2$ .
\end{enumerate}

It remains to check the consistency of all of this and reduce it to (\ref{GC}). In fact, the only thing we need to show is that $$0<M_1<M_2<\frac{2}{3}R $$ whenever $0<\Lambda<\frac{1}{12Q^2}$ i.e. $0<\Delta<1$, and $Q\ne 0$. Consider the polynomial $A(x)=x-2\Lambda x^3$. We have
$$\lim_{x\rightarrow \pm \infty}A(x)=\mp \infty$$ when $\Lambda>0$, and the roots of  $A$ are zero and $\pm a$ where $$a=\frac{1}{\sqrt{2\Lambda}}\; .$$ Also, $A$ is positive on $]0,a[$ with $R$ its local positive maximum on $x\ge 0$  and $A(R)=\frac{2}{3}R\; .$ Moreover,
\begin{align*}
0<m_1&<R<a \; ,\\
0<m_2&<R\sqrt{2}<a \; ,\\
\end{align*}
and since $A(m_i)=M_i$, if follows that $$0<M_1,M_2<\frac{2}{3}R \;.$$ Finally, to see that $M_1<M_2$ we note that
$$M_2-M_1=(m_2-m_1)(1-2\Lambda(m_2^2+m_1m_2+m_1^2))=(m_2-m_1)\left(1-\frac{2-\sqrt{1-\Delta}}{3}\right)>0\; .$$

\subsection{Photon Sphere}\label{sec:photonsphere}

Henceforth and unless otherwise specified, we will always assume that the conditions in (\ref{GC}) hold.   We will denote the three positive zeros of $f$ by $0<r_1<r_2<r_3$. The hypersurfaces $\{r=r_i\}$ for $i=1,2,3$ are respectively the {inner horizon}, the {outer horizon}, and the {cosmological horizon}. %The \emph{exterior static region}\index{Exterior static region} in which we are interested is $\mathcal{N}=\R_t \times]r_2,r_3[\times\mathcal{S}^2_{\omega}$ where $\omega=(\theta,\varphi)$.

Let us now precise what we mean by a photon sphere and continue the proof of Proposition \ref{1photonsphere}, namely, that there is only one photon sphere and it is situated in the exterior static region $\mathcal{N}=\R_t \times]r_2,r_3[ \times \mathcal{S}^2_{\omega}$. We recall the definition of the Christoffel symbols of the Levi-Civita connection determined by the RNdS metric: $$\Gamma^k_{ij}=\hf g^{kl}\left(\dl_i g_{jl}+\dl_j g_{il}-\dl_l g_{ij} \right)\; .$$
In the coordinates $(t,r,\theta,\varphi)=(x^0,x^1,x^2,x^3)$, the non zero Christoffel symbols are:
\begin{gather}
\Gamma^0_{01}=-\Gamma^1_{11}=\frac{f'}{2f} ~~~ ;~~~ \Gamma^1_{00}=\frac{ff'}{2} ~~~ ; ~~~ \Gamma^1_{22}=-rf ~~~ ; ~~~ \Gamma^1_{33}=-rf\sin(\theta)^2  \nonumber \\ \Gamma^2_{12}=\Gamma^3_{13}=\frac{1}{r} ~~~ ; ~~~ \Gamma^2_{33}=-\cos(\theta)\sin(\theta) ~~~;~~~ \Gamma^3_{23}= \cot(\theta) ~~ . \label{Crstflsymb}
\end{gather}

If we take a non zero purely rotational vector field along the angle $\varphi$ it will be of the form
$$\mathcal{V}=\alpha \frac{\partial}{\partial t} + \beta \frac{\partial}{\partial \varphi} \, ,$$
and the condition for it to be null is $g(\mathcal{V},\mathcal{V})=0$, which is

\begin{equation}\label{nullsphere}
f(r)=\frac{\beta^2 r^2 \sin(\theta)^2}{\alpha^2} \, .
\end{equation}

\begin{figure}
	\centering
	\begin{subfigure}{1\textwidth}
		\centering
		\includegraphics[width=\linewidth]{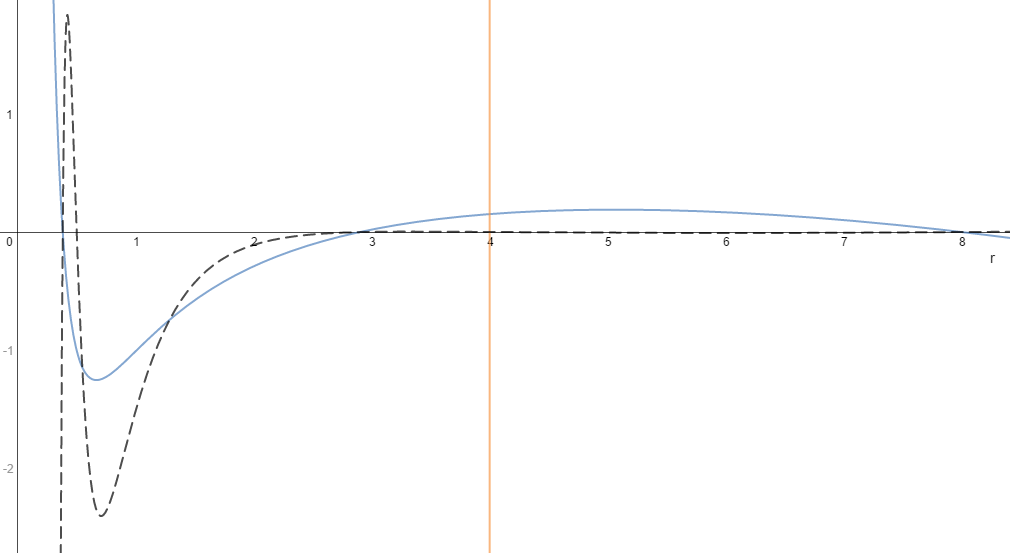}
		\caption{}
		\label{fig:fandacceleration-sub1}
	\end{subfigure}
	\begin{subfigure}{1\textwidth}
		\centering
		\includegraphics[width=\linewidth]{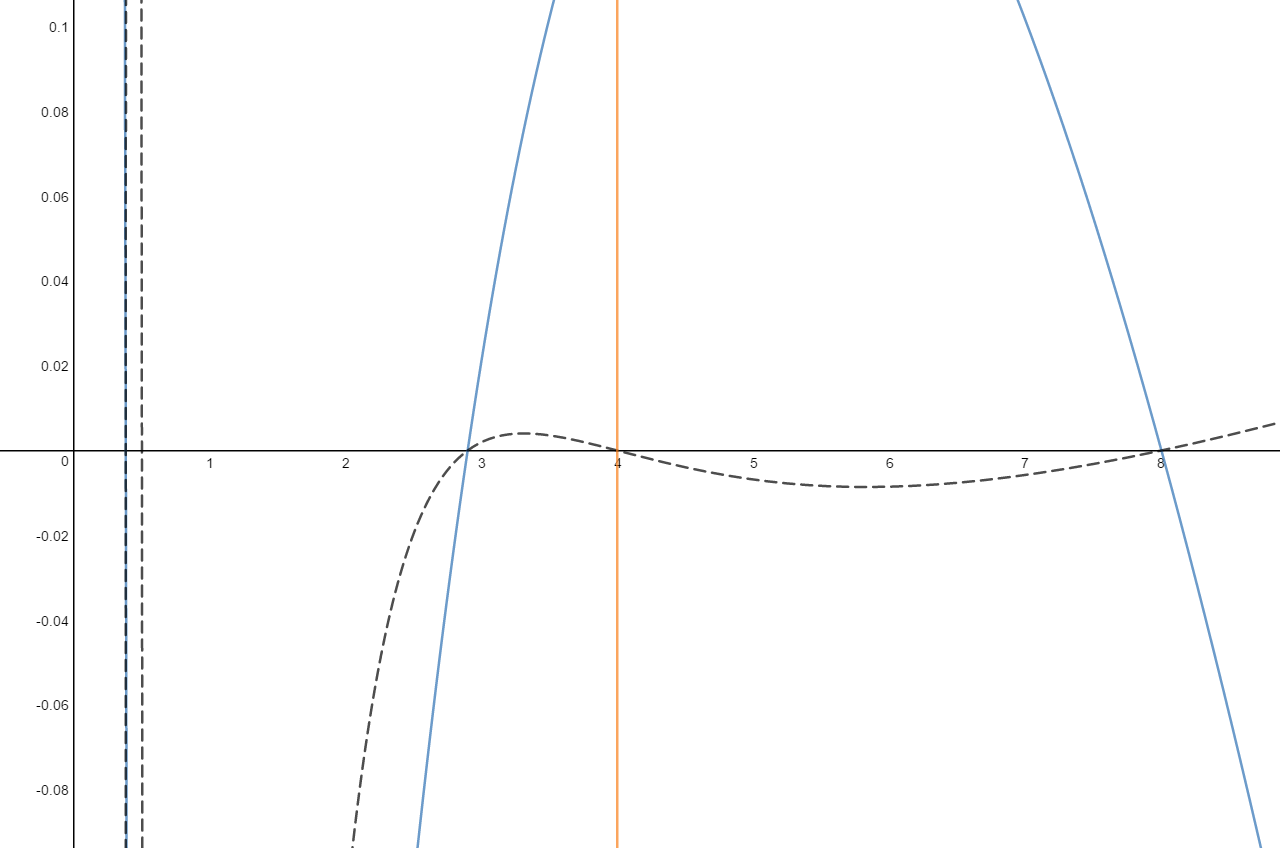}
		\caption{Rescaling the vertical axis.}
		\label{fig:fandacceleration-sub2}
	\end{subfigure}
	\caption{\emph{Numerical example: $Q=1,\; M=1.5, \; \Lambda=0.01 $. The function $f$ is the continuous curve, while the coefficient of the radial acceleration $f(2^{-1}f' - r^{-1}f)$ is the doted curve. The photon sphere is at the vertical line (r=4) $r=P_2$.}}
	\label{fig:fandacceleration}
\end{figure}
\begin{figure}
	\centering
	\includegraphics[scale=1]{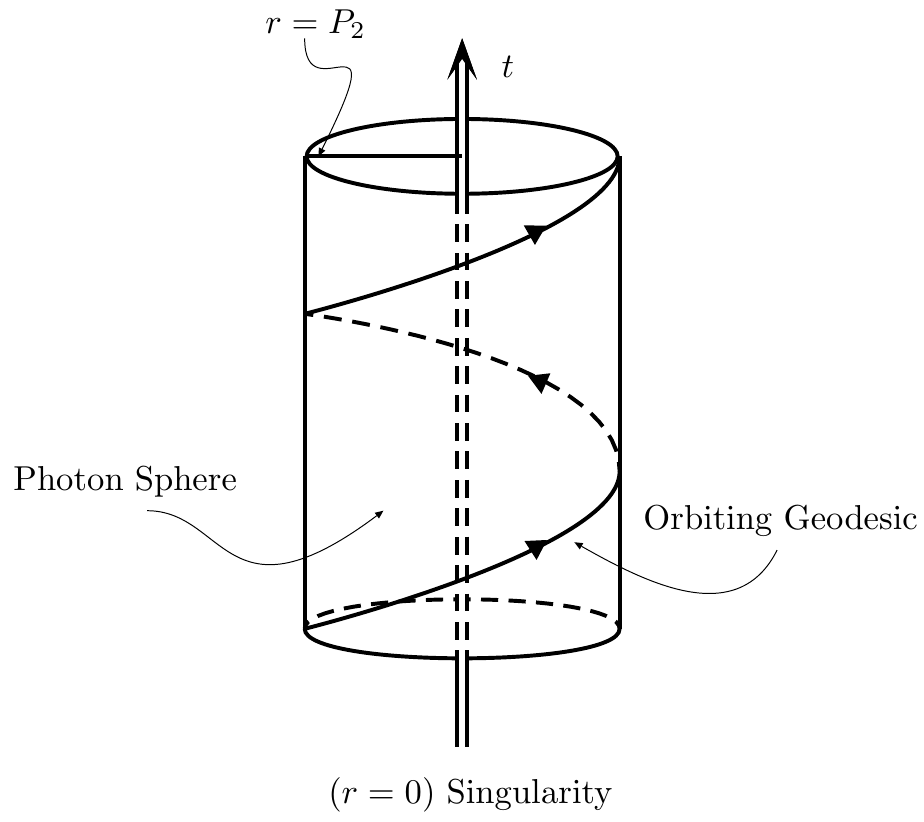}
	\caption{\emph{A null geodesic orbiting the black hole at the photon sphere $\{r=P_2\}$ which is a timelike hypersurface.}}
	\label{fig:phtnsphr}
\end{figure}

Therefore, a photon sphere could only exist in the regions where $f\geq0$, but we work in the static region outside the black hole, that is $r_3>r>r_2$, thus $f>0$. Also, condition (\ref{nullsphere}) implies that
\begin{equation}\label{rotationalvector}
\mathcal{V}=\alpha \left( \frac{\partial}{\partial t} \pm \frac{ \sqrt{f}}{r \sin(\theta)}\frac{\partial}{\partial \varphi} \right) \, ,
\end{equation}
and so it is enough to examine the case where $\alpha=1$.

Given $\mathcal{V}$ in this form, i.e. $$\mathcal{V}=\left(1,0,0,\pm \frac{ \sqrt{f}}{r \sin(\theta)}\right) \; ,$$ we calculate
\begin{eqnarray*}
	\nabla _\mathcal{V} \mathcal{V} &=& \mathcal{V}^a\left(\partial_a \mathcal{V}^c + \mathcal{V}^b \Gamma^c_{ab} \right)\partial_c \\
	&=&\Gamma^c_{00}\partial_c +\left(\mathcal{V}^3\right)^2  \Gamma^c_{33}\partial_c  \\
	&=&\Gamma^1_{00}\partial_r + \left(\mathcal{V}^3\right)^2 \left(\Gamma^1_{33}\partial_r + \Gamma^2_{33}\partial_{\theta}\right)\\
	&=& f\left(\frac{f'}{2} -\frac{f}{r}\right)\partial_r - \frac{\cot(\theta)f}{r^2} \partial_{\theta} \; .\\
\end{eqnarray*}
Thus, and since we have $f>0$ , we see that
\begin{eqnarray}\label{accelzero}
\nabla _\mathcal{V} \mathcal{V}=0 &\Leftrightarrow& f\left(\frac{f'}{2} -\frac{f}{r}\right)=0  ~~ \textrm{ and }~~ \frac{\cot(\theta)f}{r^2}=0 \nonumber\\
&\Leftrightarrow&  rf'- 2f=0 ~~ \textrm{ and }~~ \theta=\frac{\pi}{2}
\end{eqnarray}

One then can see that if we assume from the beginning that $\theta=\frac{\pi}{2}$, the integral curves of $\mathcal{V}$ at the zeros of $rf'(r)-2f(r)$ are geodesics. Hence, by the spherical symmetry, we get a full ``sphere'' of null geodesics outside the black hole, this is referred to as the {photon sphere} around the black hole.

As we can see $$rf'(r)-2f(r)=\frac{6M}{r} - \frac{4Q^2}{r^2}-2 \; ,$$ thus by studying the polynomial $$ S(r)=-r^2+3Mr-2Q^2$$ we can determine the zeros.
The discriminant of $S$ is
$$\Delta_S=9M^2-8Q^2=(3M-2\sqrt{2}Q)(3M+2\sqrt{2}Q)$$
which is positive if
$$M>\frac{2\sqrt{2}|Q|}{3} \, .$$

The two roots, if they exist, have the expressions
\begin{equation}\label{PhotonSpheres}
P_1=\frac{3M-\sqrt{\Delta_S}}{2} ~~~~ \mathrm{and} ~~~~ P_2=\frac{3M+\sqrt{\Delta_S}}{2}\; .
\end{equation}

We can directly show that the last inequality holds when (\ref{GC}) is satisfied, however, by studying the sign of $rf'-2f$ near the zeros of $f$, not only can one show that it has two zeros but also one can know their positions relative to the horizons, which is the important thing. This is Proposition \ref{1photonsphere} and the argument is in its proof which we will present now.

\begin{proof}[Continuation of the Proof of Proposition \ref{1photonsphere}]\label{proofOnephotonsphere}
	We showed that $f$ has three positive zeros $r_1,r_2,$ and $r_3$ if (\ref{GC}) holds. Note that $f$ and $P$, given in (\ref{P(r)}), are both smooth and have the same sign over $]0,+\infty[$, and we know the sign of $P$ everywhere. In a small interval around $r_1$, $f$ is decreasing since it is positive to the left of $r_1$ and negative to its right, thus $f'<0$ over this interval. Shrinking the interval if necessary, it follows that in the acceleration of the vector field $\mathcal{V}$ (see (\ref{accelzero}))  $$\nabla _\mathcal{V} \mathcal{V}=f\left(\frac{f'}{2} -\frac{f}{r}\right)\partial_r$$ the factor $f(2^{-1}f' - r^{-1}f)$ is negative to the left of $r_1$ and positive to its right. Using exactly the same logic, the last statement holds true for $r_2$ and $r_3$ also. (Figure \ref{fig:fandacceleration})
	
	Since the acceleration vector field is continuous, it must vanish in order to change sign. And since its zeros are $\{r_1,r_2,r_3,P_1,P_2\}$ (see (\ref{PhotonSpheres})), then by the above argument the zeros are necessarily ordered as follows: $r_1<P_1 <r_2<P_2<r_3$, which is what we wanted to prove. 
\end{proof}

Note that $\{r=P_1\}$ is \emph{not} a photon sphere since $f$ is negative on $]r_1,r_2[$ and so the rotational vector $\mathcal{V}$ is necessarily spacelike. This means that there are no orbits inside the black hole horizon, which is consistent with the fact that this region is dynamic. We also note that in spite the covering of the photon sphere by null geodesics it is not a null hypersurface, as a matter of fact, the spacelike  vector $\dl_r$ is a normal to the photon sphere hypersurface, and therefore it is timelike hypersurface. (Figure \ref{fig:phtnsphr})

\section{Maximal Analytic Extension}\label{sec:maximalextension}

In the RNdS coordinates $(t,r,\omega)\in\mathcal{M}=\R_t \times ]0,+\infty[_r \times \mathcal{S}^2_{\omega}$ the metric
\begin{equation}\label{metricgRNdS}
g=f(r)\d t^2-\frac{1}{f(r)}\texttt{d}r^2-r^2\d \omega^2,
\end{equation}
appears to be singular at $r=r_i$ where the factor $f^{-1}(r)$ blows up. So, the metric $g$ in these coordinates is actually defined on $\mathcal{M}$ with these hypersurfaces removed. The removal of these hypersurfaces disconnects the spacetime and divides it into four open regions $\mathrm{U}_i$, $i$ from 1 to 4. Usually, a spacetime is defined to be a connected smooth Lorentzian 4-manifold, so, we consider each of these open regions separately, and we write $g_i$ for $g|_{\mathrm{U}_i}$ when necessary. To understand the meaning of these coordinate singularities, we shall extend each of these regions by analytic extensions covering the horizons. We say that a connected {analytic Lorentzian 4-manifold}\footnote{An analytic $n$-manifold is a topological space with an atlas whose charts are analytically related, i.e. the transition maps between the charts are analytic functions of $\R^n$. Also, if it is Lorentzian, we require the metric defined on it to be analytic.} $(\tilde{\mathrm{U}},\tilde{g})$ is an {analytic extension} of $(\mathrm{U}_i,g_i)$ if the latter is isometrically embedded in the first and its image is a proper subset. Since we only consider analytic extensions, we shall refer to them simply as extensions, and an {inextendible spacetime} will be a spacetime which has no extensions.    

Intuitively, the places to extend are near the the removed hypersurfaces, and just as in other spacetimes, like Schwarzschild, there is a known way of doing it by simple changes of coordinates on each region which when extended to their maximal natural domains of definition produce extensions of the whole spacetime $\mathcal{M}$ including the hypersurfaces $\{r=r_i\}$. These are the Eddington-Finkelstein coordinates and extensions, and they are covered with single coordinate charts\footnote{Except for the fact that $\mathcal{S}^2$ needs multiple charts to fully cover it.}. However, as shall be seen, each of these extensions separately is not maximal, that is, has an extension itself. From the possible time orientations on $\mathcal{M}$ , we shall see that the different Eddington-Finkelstein extensions are in some sense complementary, and can be done all at the same time to give a bigger extension. This extension consists of an infinite number of different overlapping Eddington-Finkelstein coordinate charts, yet, there will be in the maximal extension isolated spheres of radii $r_i$s, called the {bifurcation spheres} or {crossing spheres, which are not covered by any of the Eddington-Finkelstein coordinate charts. To cover these, we need to introduce new coordinate charts on $\mathcal{M}$ and extend them, they are called the Kruskal-Szekeres coordinates. It turns out that the maximal extension can be completely covered using three families of Kruskal-Szekeres coordinate domains. 

Moreover, the maximal analytic extension of RNdS satisfies a rather stronger inextendiblity property. Besides being inextendible in the sense we described above, it is also {locally inextendible}: It has no open non-empty subset whose closure is non-compact and can be embedded in an analytic manifold with a relatively compact image. With this taken into account, there is a unique maximal connected analytic extension of RNdS which is locally inextendible, as long as we do not make identifications that change the topology.

The two dimensional diagrams presented in this section are two dimensional cross-sections of the spacetime at fixed generic angular direction $\omega_0=(\theta_0,\varphi_0)$, or equivalently, they are quotients of the spacetime by the action of the rotation group.  

\subsection{Reissner-Nordstr{\o}m-de Sitter Coordinates}

We start by reviewing some properties of the RNdS coordinates $(t,r,\omega)$. Consider the following open subsets of $\mathcal{M}$, which we also refer to by I, II, III, and IV, respectively:
\begin{eqnarray*}
	\mathrm{U}_1&=&\R_t \times ]0,r_1[_r \times \mathcal{S}^2_{\theta,\varphi}\; ; \\
	\mathrm{U}_2&=&\R_t \times ]r_1,r_2[_r \times \mathcal{S}^2_{\theta,\varphi}\; ; \\
	\mathrm{U}_3&=&\R_t \times ]r_2,r_3[_r \times \mathcal{S}^2_{\theta,\varphi}\; ; \\
	\mathrm{U}_4&=&\R_t \times ]r_3,+\infty[_r \times \mathcal{S}^2_{\theta,\varphi}\; ,
\end{eqnarray*}
and let $I_i$ be the corresponding interval of $r$ such that
\begin{equation}\label{intervalsofrI_i}
\mathrm{U}_i=\R_t \times I_i \times \mathcal{S}^2_{\theta,\varphi}\; .
\end{equation}

We orient $\mathcal{M}$ by requiring $(\dl_t, \dl_r,\dl_\theta,\dl_\varphi)$ to be a positively oriented frame. On the other hand, because $\mathcal{M}$ is not connected as a Lorentzian manifold when we remove the hypersurfaces at $r=r_i$, there is no canonical way of defining a continuous time-orientation on it a priori. For example, while $\dl_t$ is timelike in I and III, it is spacelike in II and IV where $\dl_r$ is timelike.
%When we do decay in chapter \ref{decaychapter} we choose $\dl_t$ to be future-oriented in III, but this says nothing about the time-orientation of $\dl_r$ or any other causal vector in regions II and IV. 
In other words, $\mathcal{M}$ admits more than two time-orientations. In effect, each connected component has exactly two time-orientations, $\pm \dl_t$ for I and III, and $\pm \dl_r$ for II and IV. This amounts to a total of sixteen different configurations for time-orienting $\mathcal{M}$. We shall see that each configuration is isometrically embedded, via a time-orientation preserving embedding\footnote{We say that an isometric embedding preserves time-orientation if it maps future (past) oriented causal vectors to future (past) oriented ones.}, in a connected part of the maximal extension. When we want to distinguish between different time-orientations, we shall designate $(\mathrm{U}_1,+\dl_t)$ by I, and $(\mathrm{U}_1,-\dl_t)$ by I$'$, and the same for $\mathrm{U}_3$. The time-orientation on the other regions is indicated similarly, with II and IV for $+\dl_r$, and  II$'$ and IV$'$ for $-\dl_r$.

We note that $\mathcal{M}$ admits no global timelike Killing vector field. Only regions $\mathrm{U}_1$ and $\mathrm{U}_2$ admit a timelike Killing vector field, and hence are {stationary}, and in fact, since this vector field is $\dl_t$ which is orthogonal to the foliation by the spacelike hypersurfaces $\{t=cst\}$, regions I and III are {static}. Regions II and IV are {dynamic} (not stationary) which implies, in particular, that in these regions no observer or light can ``hover'' or orbit at a fixed distance from the singularity at $r=0$. 

\begin{figure}
	\centering
	\includegraphics[scale=1]{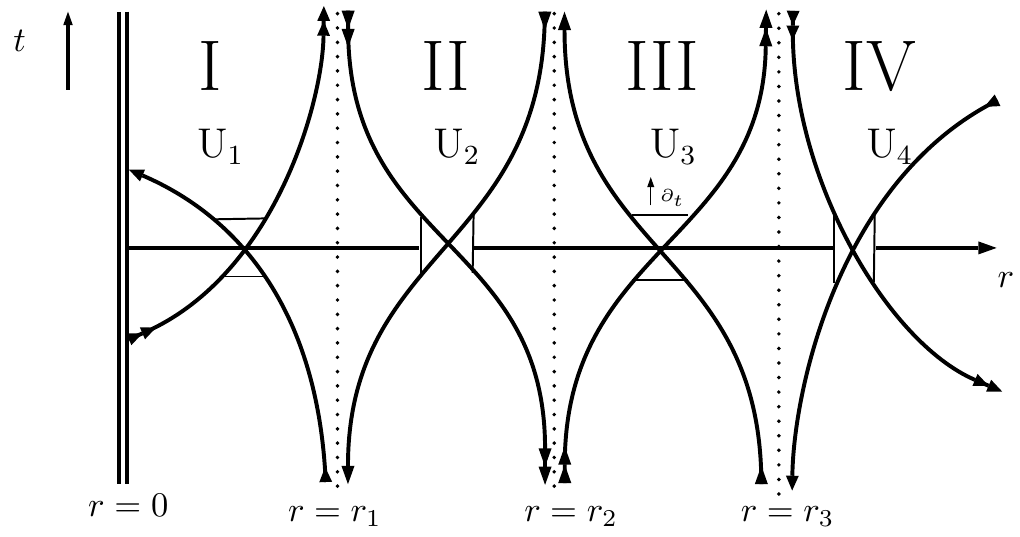}
	\caption{\emph{The open sets $\mathrm{U}_i$'s of $\mathcal{M}$, also referred to as regions} I, II, III, \emph{and} IV. \emph{ The radial null geodesics shown are integral curves of the null vector fields $Y^\mp=f^{-1}\dl_t \pm \dl_r$, with the light cones shown where they meet, and the arrows on the geodesics show the increasing direction of their affine parameters.}}
	\label{fig:RNdS}
\end{figure}

Since $(\mathcal{M},g)$ is a spherically symmetric spacetime, radial null geodesics are particularly important. First, consider a radial null geodesic $\gamma$ of $\mathcal{M}$, that is, satisfying the geodesic equations
$$\ddot{x}^k + \dot{x}^i \dot{x}^j \Gamma^{k}_{ij}=0\; ,$$
which reduce to
\begin{align*}
\ddot{t}&=- \dot{t} \dot{r} \frac{f'}{2f} \; , \\
\ddot{r}&=0 \; ,
\end{align*}
and $\gamma(s)=(t(s),r(s),\omega_0)$ for some fixed $\omega_0$, with $\dot{\gamma}(s)=\dot{t}(s)\dl_t|_{\gamma(s)}+ \dot{r}(s)\dl_r|_{\gamma(s)}$ a null vector, i.e.
$$\dot{t}^2 =\frac{\dot{r}^2}{f^{2}} \; .$$ 
Then, $r=c_1 s + c_2$ for some constants $c_1$ and $c_2$, so $r$ is an affine parameter, and we have $\dot{t}=\pm c_1 f^{-1}$. Therefore, $\gamma$ is an integral curve of a vector field of the form $c(f^{-1}\dl_t \pm \dl_r)$ for some non zero constant $c$, and hence it is sufficient to study the integral curves of the vector fields\footnote{The choice of the sign in the names of the vector fields $Y^-=f^{-1}\dl_t+\dl_r$ and $Y^+=f^{-1}\dl_t-\dl_r$ may seem strange or unpleasant at first, but we feel it is more natural this way since an integral curve of $Y^-$ is a curve of constant $u_-$ where $u_-=t-r_*$ is the retarded time coordinate of the Eddington-Finkelstein Extension, and that of $Y^+$ is given by a constant $u_+=t+r_*$. See section \ref{sec:EFextensions}} $Y^{\mp}=f^{-1}\dl_t \pm \dl_r$ that generates the others (figure \ref{fig:RNdS}). 

Two particular features of the spacetime can be seen from the radial null geodesics and the directions of the light cones in regions I and IV. The singularity at $r=0$ has a timelike nature and is more of a ``place'' in space, which can therefore be avoided. The end-points of the null geodesics at $r=\infty$, denoted by \index{$\scri$}$\scri$, can be understood as a smooth spacelike boundary by means of a conformal rescaling\footnote{$\scri$ is perhaps better described as an ``instant of time'' in the infinite future or past of an observer in region IV.}. The timelike singularity is of course due to the charge of the black hole or the Reissner-Nordstr{\o}m aspect of the spacetime, and at the other end, the spacelike null infinity is nothing but a manifestation of the cosmological constant, i.e. the de Sitter background of the spacetime.

%\begin{figure}
%\centering
%\includegraphics[scale=1]{dSRN}
%\caption{\emph{The opens $\mathrm{U}_i$'s of $\mathcal{M}$, also referred to as the regions} \textrm{I, II, III,} \emph{and} IV.\emph{ The ingoing null geodesics (discontinuous curves) and the outgoing (continuous curves) are shown with the directions (but not the angle) of the null cones represented by the small cones.}}
%\label{fig:UiRNdS}
%\end{figure}
\subsection{Regge-Wheeler Charts}\label{sec:RWcharts}

If $\gamma^{-}$ is an integral curve of $Y^-=f^{-1}\dl_t + \dl_r$, then $r$ is an affine parameter of $\gamma^-$, and $$\frac{\d(t\circ \gamma^-)}{\d r}(r)=\frac{1}{f(r)}.$$ Thus, $t(\gamma^-(r))$ is, up to an integration constant, nothing but the Regge-Wheeler coordinate function $r_*(r)$ which we will presently define, and $\gamma^-(r)=(r_*(r)+C,r,\omega_0)$ for some constant $C$. Similarly, an integral curve of $Y^+=f^{-1}\dl_t - \dl_r$ is of the form $\gamma^+(s)=(C-r_*(-s),-s,\omega_0)$ defined for $s<0$. If we choose $\dl_t$ to be future-oriented on $\mathrm{U}_3$, the null vector fields $Y^\mp$ will be future-oriented on $\mathrm{U}_3$. Thus, $\gamma^+$ is by definition an incoming null geodesic since it is future-directed and $r$ is a decreasing function of its affine parameter $s$, while $\gamma^-$ is an outgoing null geodesic for similar reasons.

The {Regge-Wheeler radial coordinate} function (also known as the {tortoise coordinate}), is defined by
\begin{equation}\label{r_*}
\frac{\d r}{\d r_*}=f(r)  ~~ \textrm{ and }~~ r_*=0 \textrm{ when } r=P_2 ~ ,
\end{equation}
where $P_2$ is the localization of the photon sphere outside the black hole given by (\ref{PhotonSpheres})\footnote{This choice of the origin for $r_*$ is convenient when we deal with decay in \cite{mokdad_decay_2016}.}. To get the explicit expression of $r_*$ in terms of $r$, let the four zeros of $f$ be $r_i$ with $r_0<0<r_1<r_2<r_3$, and let us write $f$ as 
\begin{equation*}
f(r)=-\frac{\Lambda}{r^2}\prod_{i=0}^3 (r-r_i)\; .
\end{equation*}
We integrate,
\begin{equation}\label{TortoiseCoordinate}
r_*(r)=\int_{P_2}^r \frac{1}{f(s)}\d s= \int_{P_2}^r{-\frac{s^2}{\Lambda}\prod_{i=0}^3 \frac{1}{s-r_i}\d s}=\sum_{i=0}^3 a_i \ln |r-r_i| + a \; , 
\end{equation}
where
$$a_i=-\frac{r_i^2}{\Lambda}\prod_{j\ne i}\frac{1}{(r_i-r_j)} \quad ; \; a=-\sum_{i=0}^3 a_i \ln |P_2-r_i|  \; .$$

We note that $a_0 , a_2 >0$ and $a_1,a_3<0$, $f'(r_i)=\frac{1}{a_i}$, and $\d r=f \d r_*$ on $I_i$ defined in (\ref{intervalsofrI_i}). Since $f$ has a constant sign on each interval $I_i$, each ${r_*}_i(r):={r_*}_{|_{I_i}}(r)$ is a monotonic function on $I_i$, and in fact, analytic. Thus, on each $\mathrm{U}_i$, we define the Regge-Wheeler coordinates $(t,{r_*}_i,\omega)\in \mathrm{U}_i^* = \R \times I_i^* \times \mathcal{S}^2$, where $I_i^*=r_*(I_i)$. The metric in these coordinates is
\begin{equation*}
g=f(r)(\d t^2-\d {r_*}_i^2)-r^2\d \omega^2 \; ,
\end{equation*}
where $r=r({r_*}_i)$ is the inverse function of ${r_*}_i(r)$. 

We shall usually drop the $i$ in ${r_*}_i$ (and in other coordinates later) for clarity. The ordered basis $(\dl_t,\dl_{r_*},\dl_\theta,\dl_\varphi)$ is positively oriented on $\mathrm{U}_1^*$ and $\mathrm{U}_3^*$, and negatively oriented on the other two domains. To determine the intervals $I_i^*$ we calculate the limits of ${r_*}(r)$ at the singularity $r=0$, at the horizons $r=r_i$, and at infinity $r=+\infty$.

First,
$$\lim_{r \rightarrow 0} {r_*}(r)={r_*}(0)=b\in \R ,$$
and from the signs of the coefficients $a_i$s, we have the two sided limits:
\begin{eqnarray}\label{limitsofr_*}
\lim_{r\rightarrow r_1}{r_*}(r)&=&+\infty ,\nonumber\\
\lim_{r\rightarrow r_2}{r_*}(r)&=&-\infty ,\\
\lim_{r\rightarrow r_3}{r_*}(r)&=&+\infty ,\nonumber
\end{eqnarray}
and, 
$$\lim_{r\rightarrow +\infty}{r_*}(r)=a\; ,$$ 
since as $r$ tends to $+\infty$ we have
$$r_*(r)-a\sim \sum_{i=0}^3 a_i \ln (r)=\ln (r) \sum_{i=0}^3 a_i \; ,$$
but for any four distinct non zero numbers $x_1,x_2,x_3,x_4$ we have
\begin{equation*}
\sum_{j=1}^4 x_j^2 \prod_{k \ne j}\frac{1}{(x_j-x_k)}=0 \; ,
\end{equation*}
and so, $\sum_{i=0}^3 a_i =0$. Therefore, $ I_1^* = ]b,+\infty [ $, $I_2^* =I_3^* =\R $, and $I_4^*=]-\infty , a[$.

\begin{figure}
	\centering
	\includegraphics[width=\textwidth]{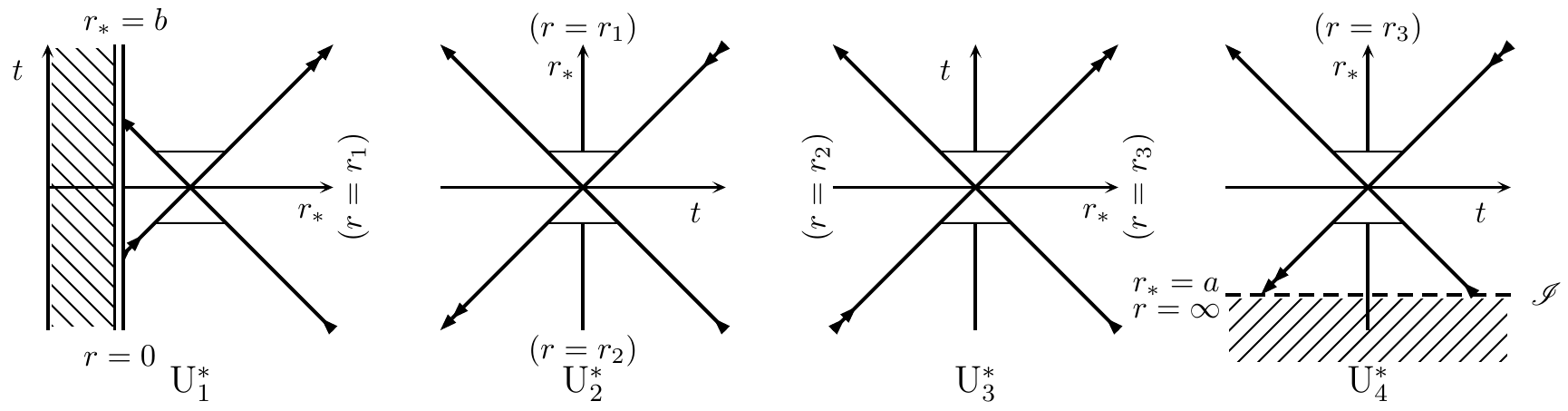}
	\caption{\emph{Oriented Regge-wheeler charts defined on $\mathrm{U}_i$. The null geodesics are integral curves of $Y^\mp=f^{-1}(\dl_t\pm\dl_{r_*})$ (lines at $\pm45$° ). The hypersurfaces $r=r_i$ (indicated in parenthesis) are off the chart since they are limits of $r_*$.}}
	\label{fig:Regge-wheelercoord}
\end{figure}

%The time-orientation of $\dl_{r_*}=f\dl_r$ is the same as $\dl_r$.
In $(t,r_*,\omega)$-coordinates,  $Y^\mp=f^{-1}(\dl_t \pm \dl_{r_*})$ and we readily see that the images of their integral curves $\gamma^\mp$ are straight lines at 45\textdegree (figure \ref{fig:Regge-wheelercoord}). The intuitive meaning of a straight line in a spacetime diagram is the worldline of a particle moving with a constant speed in the given coordinate reference frame, this is the case here too. If we consider $r_*$ as the radial coordinate, then the coordinate speeds of  $\gamma^\mp$ are indeed constants and equal to 1, i.e. $$\frac{\left(\dot{\gamma}^\mp\right)^1}{\left(\dot{\gamma}^\mp\right)^0}=1.$$ Compared to their coordinate speeds in the $(t,r,\omega)$-coordinates which is $f$, the effect of slowing down near the hypersurfaces $\{r=r_i\}$ in the $(t,r,\omega)$-coordinates is not apparent anymore in the $(t,r_*,\omega)$-coordinates but at the cost of pushing the hypersurfaces $\{r=r_i\}$ to infinite $r_*$-values. This is why $r_*$ is sometimes called the ``{tortoise coordinate}''. Finally, the boundary hypersurface $r_*=a$ is conformally spacelike, and so the spacelike nature of $\scri$ is more revealed in the Regge-Wheeler coordinates.%\footnote{In RNdS black hole, for $r$ large, $r_*$ has the opposite effect of ``bringing'' null infinity at $r=\infty$ to finite coordinate values, and it is more of ``Achilles coordinate'' there.}.

\subsection{Eddington-Finkelstein Extensions}\label{sec:EFextensions}

Let us temporarily fix a time-orientation on $\mathcal{M}$. Let $\dl_t$ be future oriented on $\mathrm{U}_1$ and $\mathrm{U}_3$, so, $ Y^\mp = f^{-1}\dl_t \pm \dl_r$ are future-oriented there, and we choose the orientation given by $\dl_r$ on $\mathrm{U}_2$ and $\mathrm{U}_4$, then $Y^-$ is future-oriented while $Y^+$ is past-oriented in regions $\mathrm{U}_2$ and $\mathrm{U}_4$. Since we defined incoming and outgoing geodesics to be future-directed, $\gamma^-$ will be an outgoing radial null geodesic on $\mathcal{M}$, while there are no incoming radial null geodesics in the dynamic regions for this particular time-orientation. In the $(t,r,\omega)$-coordinates, the coordinate expression of $\gamma^-$ has discontinuities at $r=r_i$ since its $t$-coordinate blows up because $r_*$ does, but this could be a mere bad choice of coordinate. To check this, we use the coordinates given by the flows of $Y^-$, i.e. using the geodesics $\gamma^-$ themselves as coordinate lines: For each point $p=(t_p,r_p,w_0)$ of the spacetime in the plane $\{\omega = \omega_0 \}$, there is a unique $C_p\in\R$ such that $\gamma^-_p(r)=(r_*(r) + C_p , r , \omega_0)$ passes through $p$, and $p$ can be given the new coordinates $(C_p,r_p,\omega_0)$, with $t_p=r_*(r_p)+C_p$. We thus define a new coordinate $u_-:=t-r_*$, this is Eddington-Finkelstein retarded null coordinate\footnote{A null, time, or space coordinate is one whose level surfaces are null, spacelike, or timelike hypersurfaces respectively}. The Eddington-Finkelstein retarded coordinate chart on $\mathrm{U}_3$ is 
$$({u_-}_i,r,\omega)\in \R_{{u_-}} \times {I_i}_r \times \mathcal{S}^2_{\omega}\;,$$ with ${u_-}_i=t-{r_*}_i$. In this chart the metric is: 
\begin{equation}\label{RetardedEddington-Finkelsteinmetric}
g=f(r)\d {u_-}_i^2 +2\d {u_-}_i \d r -r^2\d\omega^2\; ,
\end{equation} 
This expression of the metric is analytic for all values $(u_-,r,\omega)\in \R \times ]0,+\infty[ \times \mathcal{S}^2$, including $r=r_i$. The Lorentzian manifold $\mathcal{M}^-=\R_{u_-}\times ]0,+\infty[_r\times \mathcal{S}^2_{\omega}$\index{$\mathcal{M}^-$} with the metric (\ref{RetardedEddington-Finkelsteinmetric}) is called the {retarded Eddington-Finkelstein extension} of the RNdS manifold. Taking the orientation of $\mathcal{M}$, $(\dl_{u_-},\dl_r,\dl_\theta,\dl_\varphi)$ is positively oriented on $\mathcal{M}^-$, and when $\dl_r$ is chosen to be future-oriented\footnote{This is \emph{not} the coordinate vector field $\dl_r$ of the chart $(t,r,\omega)$. If we denote the Eddington-Finkelstein retarded coordinates by $(u_-,r_-(=r),\omega)$ then $\dl_{r_-}=f^{-1}\dl_t + \dl_r=Y^-$.}, we denote $\mathcal{M}^-$ by $\mathcal{M}^-_F$\index{$\mathcal{M}^-_F$} and call it the {future retarded extension} (figure \ref{fig:REFcoord}).  
\begin{figure}
	\centering
	\includegraphics[scale=1]{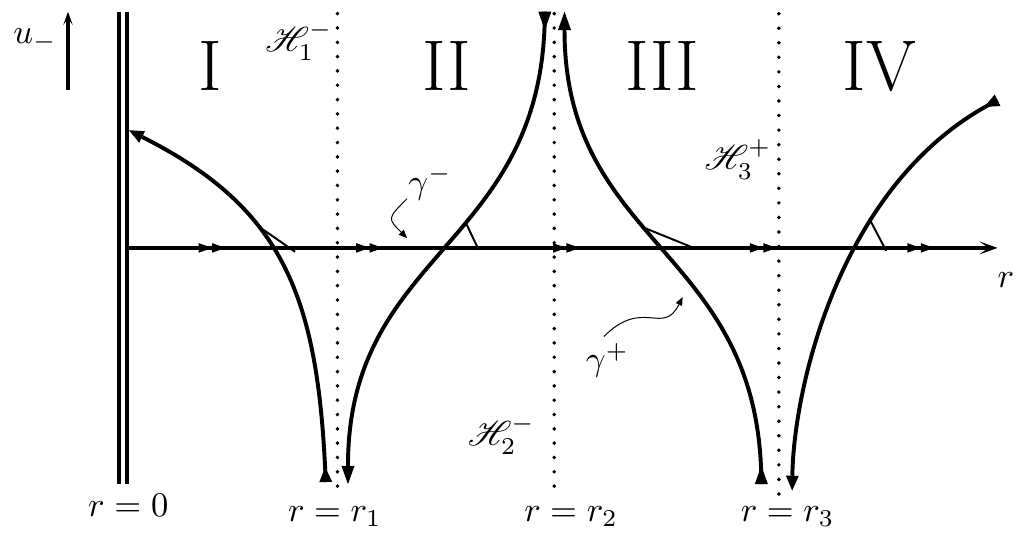}
	\caption{\emph{$\mathcal{M}^-_F$ and the integral curves of $Y^\mp$.}}
	\label{fig:REFcoord}
\end{figure}

On $\mathcal{M}^-_F$, $Y^-$ is a smooth null vector field and it is equal to  $\dl_r$ (in the retarded coordinates). Its integral curves $\gamma^-$ are outgoing radial null geodesics\footnote{Other outgoing geodesics are $-\gamma^+$ the integral curves of $-Y^+$ in the dynamic regions $\mathrm{U}_2$ and $\mathrm{U}_4$ where there are only outgoing geodesics.}, and they are just straight lines of constant $u_-$ and $\omega$: $\gamma^-(r)=(C , r , \omega_0)$. These geodesics are maximal and go through the hypersurfaces $\{r=r_i\}$ without any peculiar behaviour showing that the picture given by the RNdS coordinates is not complete and no real geometric singularities are present at $r=r_i$, and the only real singularity is at $r=0$ where the curvature becomes infinite. Since all future-directed causal curves in regions II and IV are outgoing, the hypersurfaces at $r=r_i$, the zeros of $f$, act like one way barriers which can be crossed only from the inside ($r < r_i$), and all events happening beyond them ($r > r_i$) are (for observes on the other side where $r<r_i$) hidden behind the horizons ($r=r_i$). Therefore, these hypersurfaces at $r=r_i$  are called {event horizons}, and hence $f$ is called the{horizon function} \cite{chandrasekhar_mathematical_1984}.

For an observer in III, light coming from the singularity and passing through the first two event horizons of the black hole is travelling forward in time and hence is from the past. Therefore the observer will consider the singularity to be in the past as well as the {past inner horizon} $\mathscr{H}_1^-=\R_{u_-}\times\{r=r_1\}\times \mathcal{S}^2_\omega$, and the {past outer horizon} $\mathscr{H}_2^-=\R_{u_-}\times\{r=r_2\}\times \mathcal{S}^2_\omega$, which are now regular null hypersurfaces. Similarly, the observer can only send but never receive any signal from the last horizon and $\scri$. In this extension, we denote $\scri$ by $\scri^+$ since it lies in the future of the observer, and so does the {future cosmological horizon} $\mathscr{H}_3^+=\R_{u_-}\times\{r=r_3\}\times \mathcal{S}^2_\omega$ which is a regular null hypersurface for the metric (\ref{RetardedEddington-Finkelsteinmetric}). The null horizons are generated by null geodesics each lying in a fixed angular plane (figure \ref{hrzngenerator}). This means that at the horizon some ``photons hover'' in place at $r=r_i$ and $\omega=\omega_0$. 
\begin{figure}
	\centering
	\includegraphics[scale=0.85]{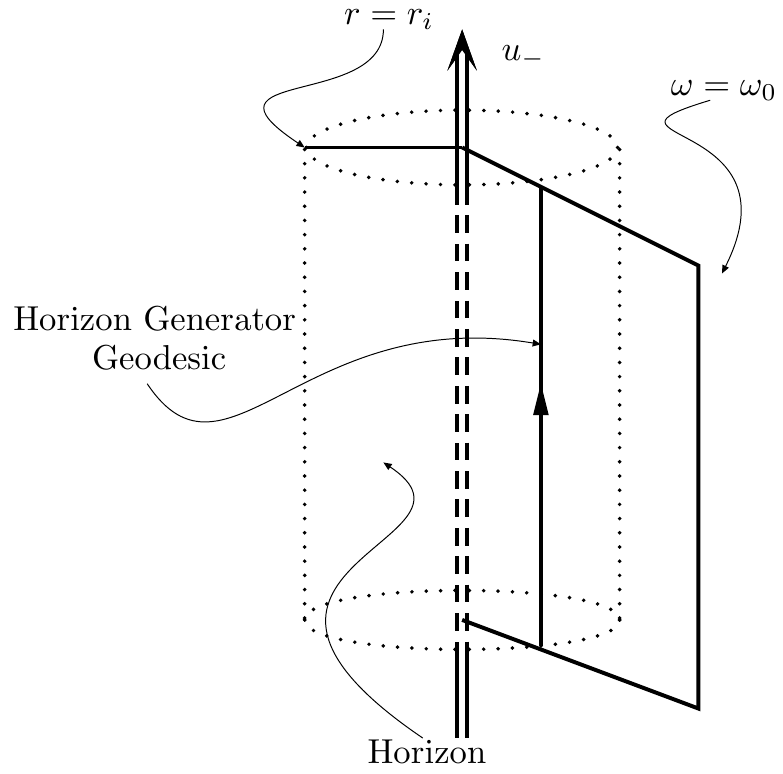}
	\caption{\emph{Integral curves of $e^{\hf f'(r_i)u_-}\dl_{u_-}$ at $r=r_i$ are null geodesics that generate the horizon $\{r=r_i\}$.}}
	\label{hrzngenerator}
\end{figure}

Although the outgoing geodesics $\gamma^-$ are inextendible in the extension $\mathcal{M}^-_F$, the incoming radial null geodesics (in the static regions I and III) $\gamma^+(s)=(C-2r_*(-s),-s+c,\omega_0)$ are not. For this reason, $\mathcal{M}^-_F$ is also called the {outgoing Eddington-Finkelstein extension}. Nonetheless, if we reverse the time-orientation, outgoing and incoming will be reversed, and the integral curves of $-Y^-$ will be the incoming geodesics crossing the horizons. We refer to $\mathcal{M}^-$ with this time-orientation as $\mathcal{M}^-_P$ the {past retarded extension}. Of course, this is not the only extension in which the incoming geodesics are inextendible, had we chosen the opposite time-orientation on $\mathrm{U}_2$ and $\mathrm{U}_4$ so that $\mathcal{M}$ is time-oriented by $\dl_t$ and $-\dl_r$ of the $(t,r,\omega)$-chart, the same procedure with $Y^+$ and $Y^-$ exchanging places would have lead instead to the {advanced Eddington-Finkelstein null coordinate} $u_+=t+r_*$ and to a new extension $\mathcal{M}^+$ covered by the single chart $(u_+,r,\omega)\in \R_{u_+} \times ]0,+\infty[_r \times \mathcal{S}^2_\omega=\mathcal{M}^+$ and endowed with the analytic metric
\begin{equation}\label{AdvancedEddington-Finkelsteinmetric}
g=f(r)\d {u_+}^2 -2\d {u_+} \d r -r^2\d\omega^2\; ,
\end{equation}
where $(\dl_{u_+},\dl_r,\dl_\theta,\dl_\varphi)$ is positively oriented and $-\dl_r$ is future-oriented. This is the {future advanced Eddington-Finkelstein extension} $\mathcal{M}^+_F$ (figure \ref{fig:AEFcoordd}). Similarly, with $\dl_r$  future-oriented we get the{past advanced Eddington-Finkelstein extension} $\mathcal{M}^+_P$. The picture given by $\mathcal{M}^+_F$ and the one given by $\mathcal{M}^-_P$ are alike but not quiet the same. In both, the singularity and the horizons at $r=r_1$ and $r=r_2$ are in the future of region III, while the horizon at $r=r_3$ is in the past of the region where also past null infinity $\scri^-$ is. In $\mathcal{M}^+_F$, we have the {future inner horizon} $\mathscr{H}_1^+=\R_{u_+}\times\{r=r_1\}\times \mathcal{S}^2_\omega$, the {future outer horizon} $\mathscr{H}_2^+=\R_{u_+}\times\{r=r_2\}\times \mathcal{S}^2_\omega$, and the {past cosmological horizon} $\mathscr{H}_3^-=\R_{u_+}\times\{r=r_3\}\times \mathcal{S}^2_\omega$. For the past extensions, $\mathcal{M}^\pm_P$, we shall denote the horizon by a minus sign when we want to be specific: $-\mathscr{H}_i^\pm$.
\begin{figure}
	\centering
	\includegraphics[scale=1]{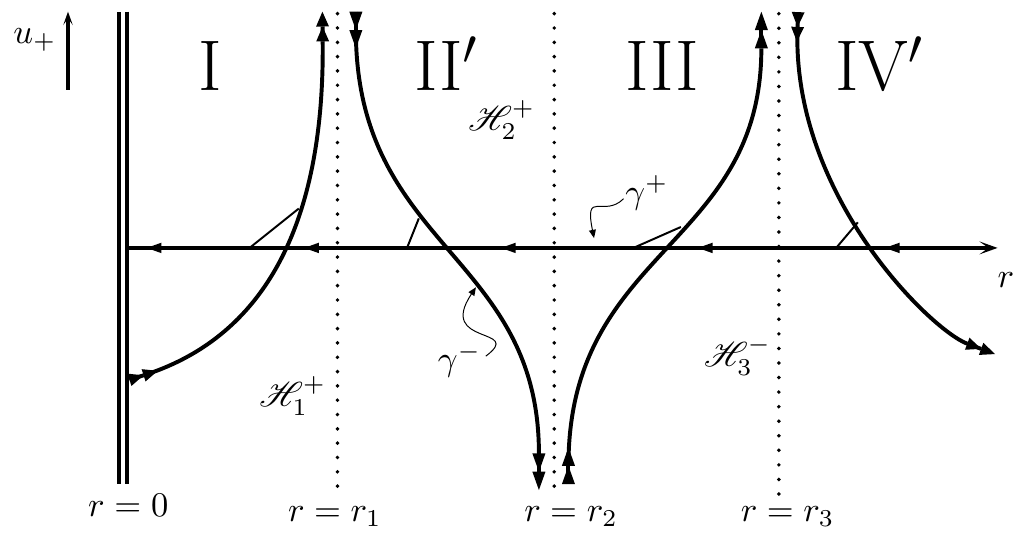}
	\caption{\emph{$\mathcal{M}^+_F$ and the integral curves of $Y^\mp$.}}
	\label{fig:AEFcoordd}
\end{figure}

With these four extensions in hand we can see what the different regions of $\mathcal{M}$ represent. Although, as seen from the geodesics $\gamma^\pm$, none of the extensions is locally inextendible and of course non is geodesically complete\footnote{In fact, even with the four combined we still do not have a geodesically complete  spacetime because of the singularity at $r=0$ beyond which the metric cannot be smoothly or even continuously extended. So we do not expect the maximal extension to be geodesically complete.}, yet, when combined they give us an almost full picture: For {almost} any radial null geodesic in $\mathcal{M}$ there is an Eddington-Finkelstein extension for which the given geodesic is future-directed and maximal i.e. extending from the singularity to $\scri$. We say almost because the null generators of the horizons are not maximal in any Eddington-Finkelstein extension, even when combined. 

From the convention of labelling the regions with primed and unprimed roman numbers to indicate time-orientation, we can easily follow the parts of different Eddington-Finkelstein extensions which are isometric in an orientation preserving manner: Different parts that carry the same label describe exactly the same geometry, the same orientation, and  the same time-orientation, so, the difference is merely a change of coordinates. Each labelled region will be covered by exactly two of these extensions, however, this is not the case for the horizons. For example, $\mathcal{M}^\pm_F$ both agree on III, but in $\mathcal{M}^-_F$ the null geodesic $\gamma^-$ intersects the past outer horizon $\mathscr{H}^-_2$ and the future cosmological horizon $\mathscr{H}^+_3$, i.e. the hypersurfaces $r=r_2$ and $r=r_3$, whereas in $\mathcal{M}^+_F$, $\gamma^-$ never touches the hypersurfaces $r=r_2$ and $r=r_3$, the future outer horizon $\mathscr{H}^+_2$ and the past cosmological horizon $\mathscr{H}^-_3$. In fact, $\mathscr{H}^+_2$ and  $\mathscr{H}^-_3$ are asymptotic to $\mathcal{M}^-_F$ where $u_-=\pm\infty$.  Thus, the horizons can not be identified with each other, this is why we distinguish between future and past horizons.

\subsection{Kruskal-Szekeres Extensions}

With the Eddington-Finkelstein extensions we {almost} have the full picture since even if we extend using the four extensions at the same time, we still do not get a locally inextendible manifold. To see why it is the case, let us examine what do we mean by doing all Eddington-Finkelstein extensions at the same time.

Each of the previous extensions is done basically by a change of coordinates on a region $\mathrm{U}_i$, then noticing that the metric in the new coordinates is analytic on a domain bigger than the original domain of the new chart, and then $\mathcal{M}$ is isometrically embedded in this bigger domain. We follow a similar strategy here, however we are not going to able cover the new extension, or even $\mathcal{M}$ for this matter, by a single coordinate chart, we need three, which is related to having three horizons. 

We start by defining on $\mathrm{U}_i$ the {double null coordinates} ${u_-}_i=t-{r_*}_i$ and ${u_+}_i=t+{r_*}_i$. We have $({u_-}_i,{u_+}_i,\omega)\in\hat{\mathrm{U}}_i$ with:
\begin{eqnarray*}
	\hat{\mathrm{U}}_1&=& \{({u_-}_1,{u_+}_1)\in\R^2 \; ;\; {u_+}_1-{u_-}_1 > 2b \} \times \mathcal{S}^2_\omega\; ; \\
	\hat{\mathrm{U}}_2&=&\R_{{u_-}_2} \times \R_{{u_+}_2}\times \mathcal{S}^2_\omega\; ; \\
	\hat{\mathrm{U}}_3&=& \R_{{u_-}_3} \times \R_{{u_+}_3}\times \mathcal{S}^2_\omega\; ; \\
	\hat{\mathrm{U}}_4&=&\{({u_-}_4,{u_+}_4)\in\R^2 \; ;\; {u_+}_4-{u_-}_4 < 2a \} \times \mathcal{S}^2_\omega\; .
\end{eqnarray*}
\begin{figure}
	\centering
	\includegraphics[scale=0.85]{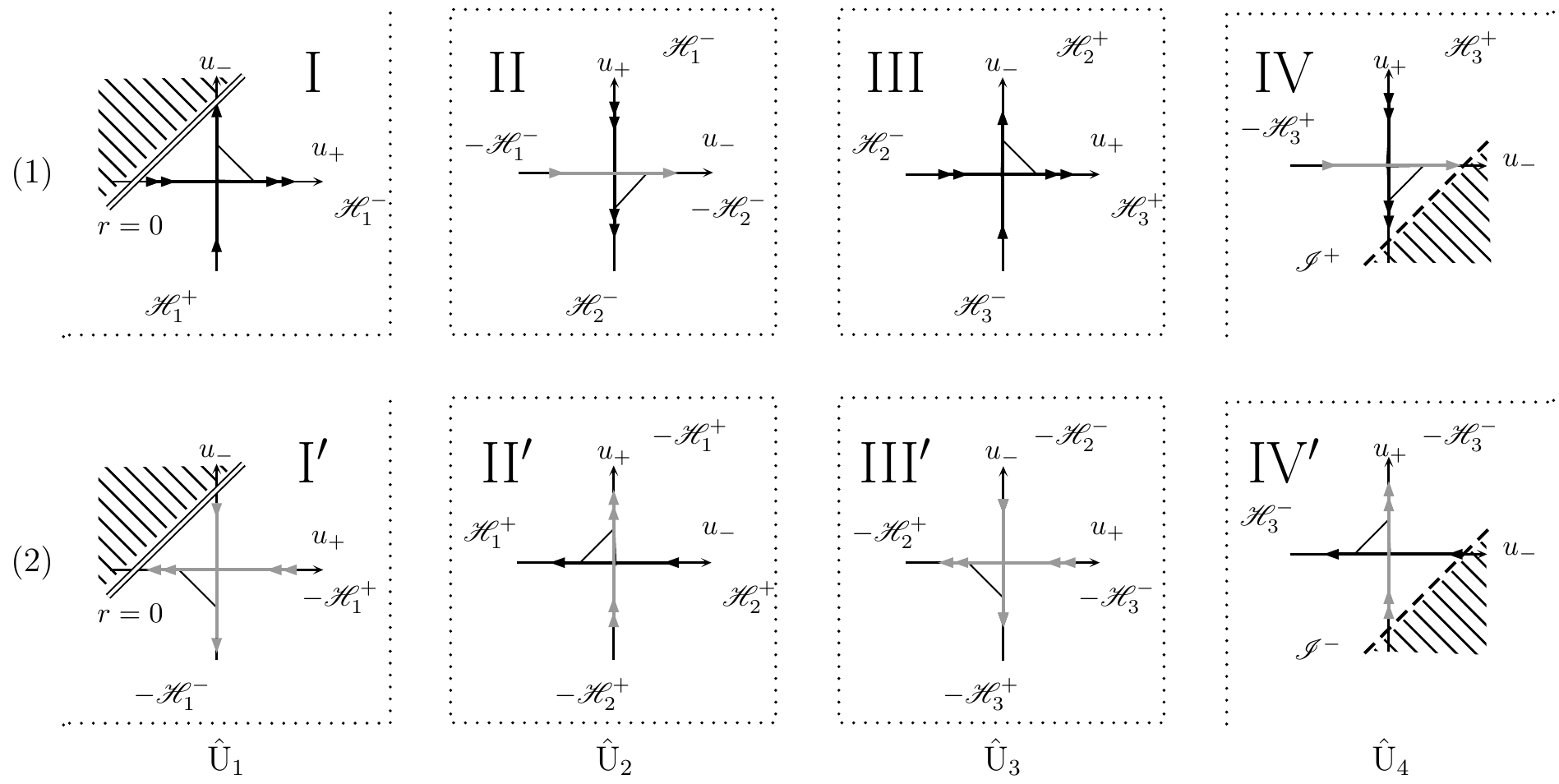}
	\caption{\emph{Oriented double null coordinates on $\hat{{\mathrm{U}}}_i$. In (1), we have time-orientation given by $\dl_t$ and $\dl_r$, while in (2) by $-\dl_t$ and $-\dl_r$. Incoming and outgoing radial null geodesics are integral curves of $Y^\mp=2f^{-1} \dl_{u_\pm}$ and of $-Y^\mp$ (shown in gray). The horizons $\pm\mathscr{H}^{\pm}_i$ (dotted lines) are asymptotic to the charts.}}
	\label{fig:doublenullcoord}
\end{figure}
The frame $(\dl_{{u_-}_i},\dl_{{u_+}_i},\dl_\theta,\dl_\varphi)$ is positively oriented on $\hat{\mathrm{U}}_1$ and $\hat{\mathrm{U}}_3$, and negatively oriented on $\hat{\mathrm{U}}_2$ and $\hat{\mathrm{U}}_4$. The metric in these coordinates is,
\begin{figure}
	\centering
	\includegraphics[width=\textwidth]{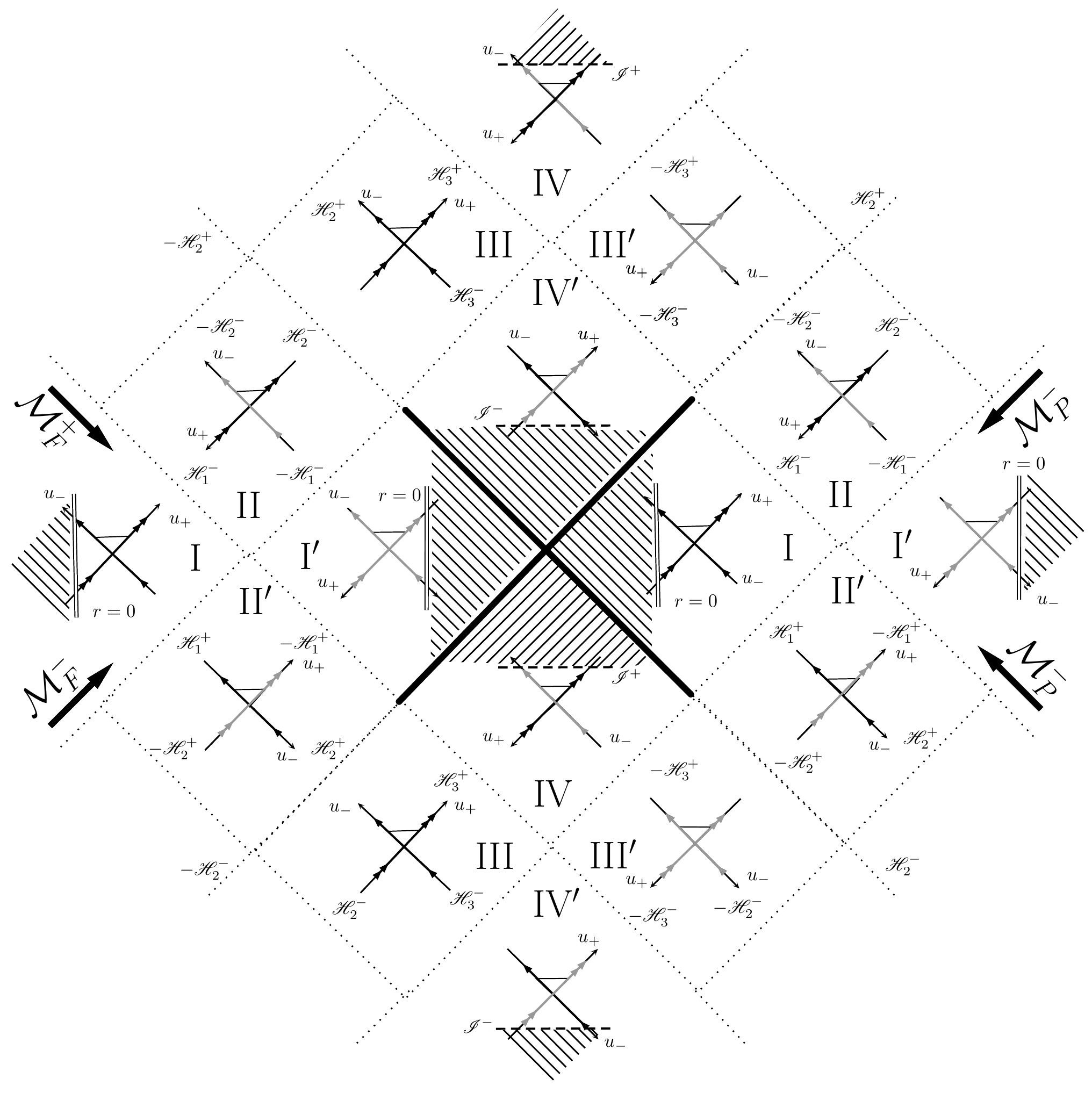}
	\caption{\emph{The gluing of double null coordinates along their asymptotic horizons, covering the different Eddington-Finkelstein extensions, with the pattern repeating infinitely. See figure \ref{fig:doublenullcoord} for the legend.}}
	\label{fig:glueingcharts}
\end{figure}
\begin{equation}\label{metricdoublenullcoord}
g=f(r)\d {u_-}_i \d {u_+}_i-r^2\d \omega^2 \; ,
\end{equation}
where $r$ is implicitly given by ${u_+}_i - {u_-}_i=2{r_*}_i(r)$.

To put these charts in the context of the Eddington-Finkelstein extensions, we need to choose orientations on the $\hat{\mathrm{U}}_i$s. We can then use the radial null geodesics and the extensions $\mathcal{M}^\pm_F$ and $\mathcal{M}^\pm_P$ to determine the asymptotics of the oriented double null coordinates charts. This is figure \ref{fig:doublenullcoord}. The next step is to ``glue'' the charts along their common (asymptotic) horizons, and since we want to understand these charts as oriented coordinate systems on Eddington-Finkelstein extensions, there is only one way of putting them together, which is shown in figure \ref{fig:glueingcharts}.

It is clear that we have left out where the four null hypersurfaces $\pm\mathscr{H}^\pm_i$ (for the same $i$) meet, this is a sphere $\mathcal{S}_i$ called the {bifurcation sphere} since the hypersurface $r=r_i$ bifurcates into four horizons. To see that the missing spheres are regular and the metric can be analytically extended on them, we need to define new coordinates on $\hat{\mathrm{U}}_i$ for which we can identify (glue) the corresponding horizons as regular hypersurfaces and not just asymptotically. On the one hand, to bring the horizons back to finite coordinate values, one choice is exponential functions of the null coordinates $u_-$ and $u_+$ with the correct weights, this is the Kruskal-Szekeres choice of coordinates. On the other hand, the metric $g_{\Omega_0}$, defined on the two dimensional space $\{\omega=\omega_0\}$ and locally given by $f(r)\d {u_-}_i \d {u_+}_i$, is locally conformally flat as we can see form the double null coordinates expression. The only coordinate transformation that leaves it in such a form, is if one of the new coordinates is function  of $u_-$ only, and the other is function of $u_+$ only. The simplest of such transformations would be
\begin{equation*}
{U_+}_i={\beta_+} e^{\alpha_+ {u_+}_i}  \quad \mathrm{and} \quad {U_-}_i=\beta_- e^{\alpha_- {u_-}_i} \; ,
\end{equation*}
with non zero constant weights $\alpha_\pm$ and $\beta_\pm$ (indexed by $i$ but we drop it for now).
%
%1-THESE ARE $\pm\alpha_i\dl_{U_\pm}$ ON $r=r_i$, SO ITS INTEGRAL CURVES $(-\alpha_i s +C_1, C_2, \omega_0)$ ARE MAXIMAL i.e. DEFINED ON $\R$ WITH $s$ AFFINE PARAMETER, CHECK IT! 
%
%2- THE MISSING SPHERES ARE INTERSECTED ONLY BY THE GEODESIC GENERATORS OF THE HORIZONS. AND THIS IS WHERE THE NEW COORDINATE IS NEEDED. 
%
%the existence of geodesics that are not maximal, i.e. their affine parameter is defined on a proper subset of $\r$, is a good indication for local extendiblity. you use them to difine a new coordinate, and then it turns out that the spacetime is extendible in a straightforward way when written in the new coordinates.
We have,
\begin{equation*}
\d {U_+}_i = \alpha_+  {U_+}_i \d {u_+}_i \quad ; \quad \d {U_-}_i = \alpha_-  {U_-}_i \d {u_-}_i \; ,
\end{equation*}
so
\begin{equation*}
\d {U_+}_i\d {U_-}_i = \alpha_+ \alpha_-  {U_+}_i {U_-}_i \d {u_+}_i  \d {u_-}_i \; .
\end{equation*}
Thus, the metric would be
\begin{equation*}
g=\frac{f(r)}{\alpha_+\alpha_- {U_+}_i  {U_-}_i}\d {U_+}_i \d {U_-}_i - r^2 \d \omega^2 \; ,
\end{equation*}
and we need to express $r$ in terms of ${U_+}_i$ and ${U_-}_i$. In fact,
\begin{equation}
{U_+}_i {U_-}_i=\beta_+\beta_- e^{(\alpha_+ - \alpha_-) {r_*}_i}e^{(\alpha_+ + \alpha_-) t},
\end{equation}
and in order to define $r$ as a function of $({U_-}_i , {U_+}_i)$ using this relation, we must eliminate the $t$ variable from the relation, so, we take $\alpha_+=-\alpha_-=:\alpha_i $. This simplifies the above expression,
\begin{equation}\label{rinU+U-}
{U_+}_i {U_-}_i=\beta_+\beta_- e^{2\alpha_i {r_*}_i}=\beta_+\beta_- A_i\prod_{j=0}^3 |r-r_j|^{2\alpha_i{a_j}}=:h_i(r),
\end{equation}
where $A_i=e^{2\alpha_i a}$. $h_i$ is a bijective function of $r$ defined on $I_i$, since ${r_*}_i$ is bijective. Thus, $r({U_-}_i , {U_+}_i)=h_i^{-1}({U_+}_i {U_-}_i)$ is a one-to-one function from $h_i(I_i)$ onto $I_i$, and actually analytic since,
$$\frac{\d h_i}{\d r}(r)= \beta_+\beta_-\frac{2\alpha_i}{f(r)}e^{2\alpha_i {r_*}_i(r)}\neq 0.$$
It follows that 
\begin{equation}\label{metricinU+U-}
g=\frac{f(r)}{-\alpha^2_i {U_+}_i  {U_-}_i}\d {U_+}_i \d {U_-}_i - r^2 \d \omega^2 \; ,
\end{equation}
is also analytic on the domain
\begin{equation*}
\mathrm{U}'_i=\{({U_-}_i ,{U_+}_i)\in \R^2; \beta_+{U_+}_i>0 ; \beta_-{U_-}_i>0  ; {U_+}_i {U_-}_i \in h_i(I_i)\}\times \mathcal{S}^2_\omega \; .
\end{equation*}
If we want $(\dl_{{U_-}_i},\dl_{{U_+}_i},\dl_\theta,\dl_\varphi)$ to be positively oriented everywhere, then we are bound to $\beta_+\beta_-<0$ on $\mathrm{U}_1$ and $\mathrm{U}_3$, and $\beta_+\beta_->0$ on $\mathrm{U}_2$ and $\mathrm{U}_4$. This means that $h_1 ~, h_3$ are negative and $h_2 ~, h_4$ are positive. There is no serious restriction in assuming that $|\beta_\pm|=1$, since it is their sign which is interesting to us. Accordingly, we have
$$h_i(r)=(-1)^i e^{2\alpha_i {r_*}_i}=(-1)^i A_i\prod_{j=0}^3 |r-r_j|^{2\alpha_i{a_j}}.$$

As before, we have defined a new coordinate system, we now try to extend its domain of definition, keeping in mind that we wish to assign finite double null coordinates for the horizons. We see from the expression of $h_i$ that for a good choice of $\alpha_i$, $h_i$ (and hence the domain of the chart) can be extended analytically to an interval containing a horizon at $r=r_j$ where $r_j$ is a boundary point of $I_i$ different than zero. Thus, the choice of $\alpha_i$ is self suggesting as $1/a_j$ for some $j$. However, if we use $1/a_j$, we shall run into trouble when extending the metric since ${U_+}_i {U_-}_i=h_i(r)$ will be zero at $r=r_j$ and $h_i$ will contain one more power of $(r-r_j)$ than $f$, so the metric will blow up. We thus take   
$\alpha_i=\frac{1}{2 a_j}.$

Therefore, for $i=1,2,3$, let $\alpha_i=\alpha_{i+1}=\frac{1}{2 a_i}$, then $A_i=A_{i+1}$ and the function
\begin{equation*}
H_i(r)=(-1)^{i}A_i(r_i-r)\prod_{j\ne i , j=0}^3 |r-r_j|^{\frac{a_j}{a_i}}=\left\{\begin{array}{cc}
h_i(r) & r\in I_i  \\
0 & r=r_i \\
h_{i+1}(r) & r\in I_{i+1}
\end{array}\right. \; ,
\end{equation*}
is continuous on $I_i\cup\{r_i\}\cup I_{i+1}$. Moreover,  since $f$ has the opposite sign of $a_i$ over $I_i$, and the same sign over $I_{i+1}$, $H_i$ is monotonic on its domain:
\begin{eqnarray*}
	\frac{\d H_i}{\d r}|_{I_i} &=& \frac{(-1)^{i}}{a_i f}e^{\frac{1}{ a_i} {r_*}_i(r)} \\
	\frac{\d H_i}{\d r}(r_i)&=& (-1)^{i+1}A_i\prod_{j\ne i , j=0}^3 |r_i-r_j|^{\frac{a_j}{a_i}}  \\
	\frac{\d H_i}{\d r}|_{I_{i+1}} &=&\frac{(-1)^{i+1}}{a_i f}e^{\frac{1}{ a_i} {r_*}_{i+1}(r)} \; ,
\end{eqnarray*}
so, $H_1$ and $H_3$ are increasing, while $H_2$ is decreasing. Thus, $H_i$ is an analytic bijection from $I_i\cup\{r_i\}\cup I_{i+1}$ onto its image, and its inverse is also analytic. To find the domain of the inverse function we take the limits. From the limits of $r_*$ in (\ref{limitsofr_*}) we have:
\begin{equation*}
-\infty<\lim_{r\rightarrow 0} H_1(r)=H_1(0)=-e^{\frac{b}{ a_1}}:=B<0 \; ,
\end{equation*}
we also have
\begin{equation*}
\lim_{r\rightarrow r_2} H_1(r)=+\infty \; .
\end{equation*}
Thus, $H_1:]0,r_2[\longrightarrow ]B,+\infty[$ . Similarly,
\begin{eqnarray*}
	\lim_{r\rightarrow r_1} H_2(r)&=&+\infty \; , \\
	\lim_{r\rightarrow r_3} H_2(r)&=&-\infty \; , \\
\end{eqnarray*}
so, $H_2:]r_1,r_3[ \longrightarrow ]-\infty,+\infty[$. Also,
\begin{eqnarray*}
	\lim_{r\rightarrow r_2} H_3(r)&=&-\infty \; ,\\
	\lim_{r\rightarrow +\infty} H_3(r)&=& e^{\frac{a}{ a_3}}=:A>0  \; ,
\end{eqnarray*}
and $H_3:]r_2,\infty[\longrightarrow ]-\infty,A[$.
\begin{figure}
	\centering
	\includegraphics[width=\textwidth]{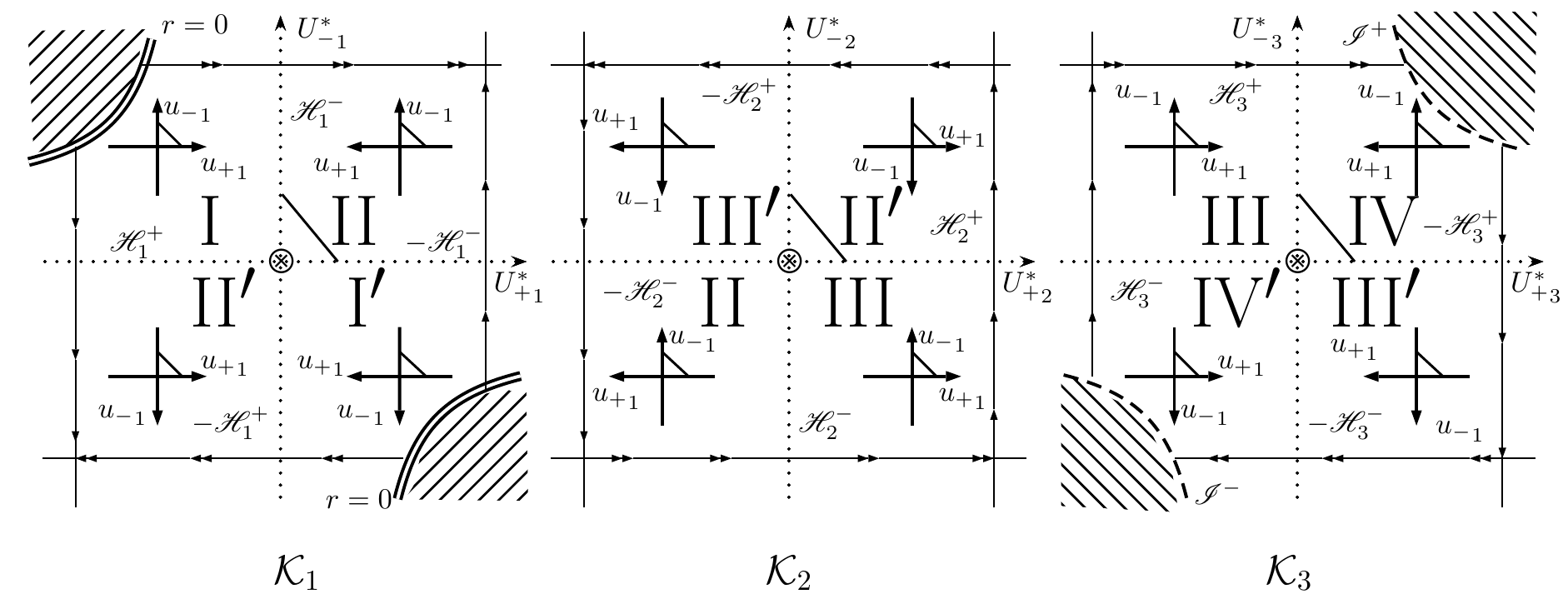}
	\caption{\emph{The Kruskal-Szekeres extensions $\mathcal{K}_i$s with the radial null geodesics $\gamma^\pm$.}}
	\label{fig:Kruskalextensions}
\end{figure}
Using the $H_i$s and the formal expression (\ref{metricinU+U-}), we can define three Lorentzian manifolds $(\mathcal{K}_i,g_{\mathcal{K}_i})$ for $i=1,2,3$, called the Kruskal-Szekeres extensions, as follows:
\begin{equation*}
\mathcal{K}_i=\{({U^*_-}_i ,{U^*_+}_i)\in \R^2 ;{U^*_+}_i {U^*_-}_i \in H_i(I_i\cup\{r_i\}\cup I_{i+1})\}\times \mathcal{S}^2_{\omega} \; ,
\end{equation*}
and
\begin{equation}\label{metricKruskal}
g_{{\mathcal{K}_i}}=\frac{-4 a_i^2 f(r)}{ H_i(r)}\d {U^*_+}_i \d {U^*_-}_i - r^2 \d \omega^2 \; ,
\end{equation}
where $r({U^*_-}_i ,{U^*_+}_i)=H_i^{-1}({U^*_+}_i {U^*_-}_i)$. Note that $g_{\mathcal{K}_i}$ is indeed analytic since the factor $(r-r_i)$ in $H_i(r)$ in the denominator is cancelled out by the same factor coming from $f(r)$, and thus the metric is regular on $r=r_i=\{{U^*_+}_i {U^*_-}_i=0\}$, in particular, it is regular on the bifurcation sphere $({U^*_-}_i, {U^*_+}_i)=(0,0)$.

To see these manifolds as local extensions of the Eddington-Finkelstein manifolds and of $\mathcal{M}$, let us embed the $\mathrm{U}_i$s in them via the transformation 
\begin{equation*}
{U^*_+}_i={\beta_+}_j e^{\frac{1}{2a_i} {u_+}_j}  \quad \mathrm{and} \quad {U^*_-}_i={\beta_-}_j e^{-\frac{1}{2a_i} {u_-}_i} \; ,
\end{equation*}
for $j=i\; ,\; i+1$, where ${u_\pm}_i=t\pm {r_*}_i$. If we want the transformation to be orientation preserving with $\mathrm{U}_i$ oriented by $(\dl_t,\dl_r,\dl_\theta,\dl_\varphi)$ which is positively oriented, then as we mentioned above, we must have ${\beta_+}_i {\beta_-}_i$ positive for $i=2,4$ and negative for $i=1,3$. Then, form the definition of $r({U^*_-}_i, {U^*_+}_i)$ and  $H_i(r)$, we see that two ``diagonally opposite quadrants'' of $\mathcal{K}_i$ are each isometric to $\mathrm{U}_i$, and the other two ``quadrants'' to $\mathrm{U}_{i+1}$, and the horizons at $r=r_i$ corresponds to the ``axis'' of $\mathcal{K}_i$, of course each of these parts of $\mathcal{K}_i$ is a product with $\mathcal{S}^2$. We note also that since $H_i(r)$ and $f(r)$ have opposite signs, $\dl_{{U^*_-}_i}+\dl_{{U^*_+}_i}$ is timelike on $\mathcal{K}_i$. The choice of this vector being future or past oriented is equivalent to fixing the sign of each ${\beta_\pm}_j$. These choices can be decided alternatively and equivalently by following the geodesics of $Y^\pm$ guided by figure \ref{fig:glueingcharts}, where $Y^\pm$ are now given in the Kruskal-Szekeres coordinates by
$$Y^\pm = \frac{1}{a_i f(r)} {U^*_\mp}_i \dl_{{U^*_\mp}_i}~.$$ 
We note that since
$$Y^- = \frac{1}{a_i f(r)} {U^*_+}_i \dl_{{U^*_+}_i}=\frac{H_i(r)}{a_i f(r){U^*_-}_i} \dl_{{U^*_+}_i},$$ 
$Y^-$ is actually defined and smooth on $\mathcal{K}_i \setminus \{{U^*_+}_i =0 \}$. Similarly for $Y^+$. The geodesics along the horizons are given by $Y^\mp_{\mathscr{H}_i}=\pm\frac{1}{2a_i}\dl_{{U^*_\pm}_i}$ on ${U^*_\mp}_i=0$. 
Figure \ref{fig:Kruskalextensions} summarizes all of this when $\dl_{{U^*_-}_i}+\dl_{{U^*_+}_i}$ is future-oriented. We remark that using this coordinates change we can recover $t$ as a function of $({U^*_-}_i ,{U^*_+}_i)$ through
\begin{equation}\label{tinU+U-}
\frac{{U^*_+}_i}{{U^*_-}_i}={\beta_+}_j{\beta_-}_j e^{\frac{t}{2 a_i}}.
\end{equation}

\subsection{The Maximal Extension}

The Kruskal-Szekeres extensions can be used as an atlas for the maximal extension. The maximal analytic extension of $\mathcal{M}$ is a Lorentzian manifold $\mathcal{M}^*$ covered by an atlas $\mathfrak{A}^*$ consisting of coordinate charts given by the $\mathcal{K}_i$s, and is endowed with the metric $g^*$ given locally as $g_{{\mathcal{K}_i}}$ (or simply $g$). The manifold $\mathcal{M}^\ast$ %which we shall consider is in fact a Penrose conformal diagram (cross-product the two sphere), in other words, it is (conformal) compactification of%
is in essence the collection of overlapping Kruskal-Szekeres extensions where the corresponding regions are identified. Equipped with the usual topology, let
\begin{equation*}
\mathcal{M}^\ast=\left(\R^2 \setminus \left(\bigcup_{k,l\in\mathbb{Z}}S_{k,l}\right)\right)\times \mathcal{S}^2 \; ,
\end{equation*}
where $S_{k,l}$ is the square block
\begin{equation*}
S_{k,l}=\left\{(x,y)\in \R^2 ; \frac{\pi}{2}\le x\sqrt{2}-2k\pi\le 3\frac{\pi}{2} \; ; \;  -\frac{\pi}{2}\le y\sqrt{2}-2l\pi\le \frac{\pi}{2}\right\} \; ,
\end{equation*}
and let the atlas be
\begin{equation*}
\mathfrak{A}^\ast=\{(A_{k,l},\phi_{k,l}),(B_{k,l},\chi_{k,l}),(C_{k,l},\psi_{k,l}) \; ; \; k,l\in\mathbb{Z}\}\; ,
\end{equation*}
with the charts defined as follows: Let $n=l-k$ and $m=l+k$, and set
\begin{equation*}
X=\frac{1}{\sqrt{2}}(y+x) \quad ; \quad Y=\frac{1}{\sqrt{2}}(y-x) \; ,
\end{equation*}
the opens\footnote{Here, we ignore the fact that the 2-sphere needs multiple charts to cover it.} $A_{k,l},B_{k,l}$, and $C_{k,l}$ are

\begin{align*}
A_{k,l}&= \left\{(x,y)\in \R^2 ; \tan(X)\tan(Y)>-1 \; ; \;  -\frac{\pi}{2}< X - m\pi< \frac{\pi}{2} \; ; \;  -\frac{\pi}{2}< Y -n\pi < \frac{\pi}{2}\right\}\times \mathcal{S}^2 \\
B_{k,l}&= \left\{(x,y)\in \R^2 ;  -\frac{\pi}{2}< X - \left(m+\hf\right)\pi< \frac{\pi}{2} \; ; \;  -\frac{\pi}{2}< Y -\left(n+\hf\right)\pi < \frac{\pi}{2}\right\}\times \mathcal{S}^2 \\
C_{k,l}&= \left\{(x,y)\in \R^2 ; \tan(X)\tan(Y)<1 \; ; \;   \right.\\ 
&\hspace{5cm}\left. -\frac{\pi}{2}< X - (m+1)\pi< \frac{\pi}{2} \; ; \;-\frac{\pi}{2}< Y -n\pi < \frac{\pi}{2}\right\}\times \mathcal{S}^2,
\end{align*}
and the chart bijections are
\begin{eqnarray*}
	\phi_{k,l}:A_{k,l}&\longrightarrow& \mathcal{K}_1 \; ; \\
	\chi_{k,l}:B_{k,l}&\longrightarrow& \mathcal{K}_2 \; ;\\
	\psi_{k,l}:C_{k,l}&\longrightarrow& \mathcal{K}_3 \; ,
\end{eqnarray*}
given by
\begin{eqnarray*}
	\phi_{k,l}(x,y,\omega)&=&({U_-}_1 ,{U_+}_1,\omega)=(\sqrt{-B}\tan(Y),\sqrt{-B}\tan(X),\omega) \; ; \\
	\chi_{k,l}(x,y,\omega)&=&({U_-}_2 ,{U_+}_2,\omega)=\left(\tan\left(Y-\frac{\pi}{2}\right),\tan\left(X-\frac{\pi}{2}\right),\omega\right) \; ;\\
	\psi_{k,l}(x,y,\omega)&=&({U_-}_3 ,{U_+}_3,\omega)=(\sqrt{A}\tan(Y),\sqrt{A}\tan(X),\omega) \; .
\end{eqnarray*}

The metric on $\mathcal{M}^\ast$ is $g^\ast$ whose coordinate expression on each chart domain is given by (\ref{metricKruskal}). This extension is the maximal analytic extension of RNdS manifold. It is maximal in the sense that it is locally inextensible. It is also unique if the topology is not changed.

\begin{figure}
	\hspace{-0.5cm}
	\includegraphics[width=\textwidth]{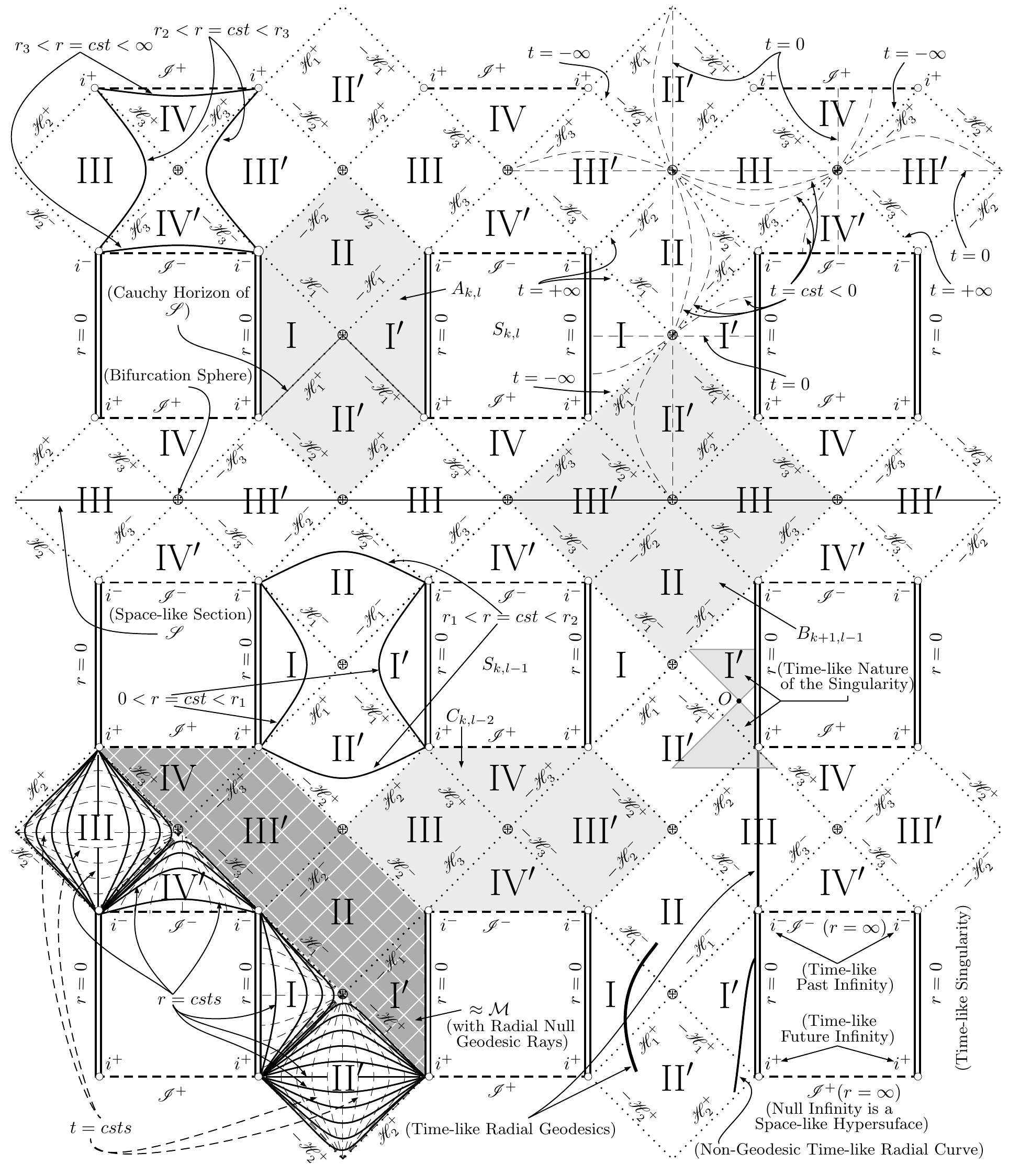}
	\caption{\emph{The Structure of $\mathcal{M}^\ast$.}}
	\label{fig:OpensConnected1}
\end{figure}

The structure of $\mathcal{M}^\ast$ is shown in figure \ref{fig:OpensConnected1}. First, we note that the metric is analytic and well behaved at all points of $\mathcal{M}^\ast$, including the horizons which are now given by ${U_+}_i {U_-}_i=0$ and the bifurcation spheres $({U_+}_i ,{U_-}_i)=(0,0)$. The RNdS radius $r$ is a scalar field on $\mathcal{M}^*$, but the same can not be said regarding the time parameter $t$, which is given in the different regions through (\ref{tinU+U-}) as shown in the diagram, and is not defined on the horizons where it becomes infinite. 

$\mathcal{M}^\ast$ contains infinitely many isometric copies of the original spacetime $\mathcal{M}$. Each consists of four regions numbered from one to four in roman, possibly primed or mixed. There are sixteen (infinite) families of these copies, each family corresponds to one of the sixteen different ways of time-orienting $\mathcal{M}$. Four of the families correspond to the Eddington-Finkelstein extensions of $\mathcal{M}$. Examples of the others along with these four are shown in figure \ref{fig:timeorientations16}. 

The causal structure of $\mathcal{M}^\ast$ can also be seen from figure \ref{fig:glueingcharts}.
\begin{figure}
	\centering
	\includegraphics[scale=0.7]{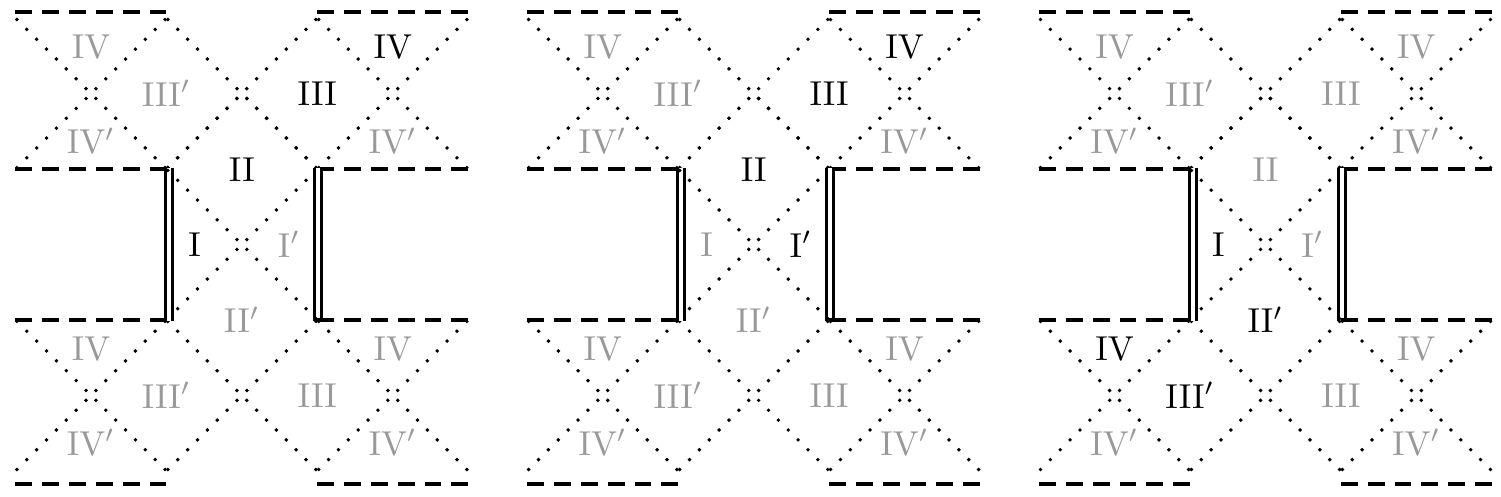}
	\caption{\emph{Examples of different time-orientations on $\mathcal{M}$ as connected subsets of $\mathcal{M}^*$.}}
	\label{fig:timeorientations16}
\end{figure}
Upon choosing a time orientation on $\mathcal{M}^\ast$, say $\dl_{{U^*_-}_i}+\dl_{{U^*_+}_i}$, then all future directed timelike causal curves in region IV end at $\scri^+$, and all past directed causal curves end at $\scri^-$. Unlike Minkowski, Schwarzschild, Reissner-Nordstr{\o}m, or Kerr spacetimes, in RNdS, null infinity or $\scri$ is not a null ``hypersurface'', instead it is spacelike due to the De Sitter nature of our spacetime.  Using the conformal factor $\sqrt{|f^{-1}|}$ one can define the metric on this hypersurface, and see that it is indeed spacelike for the conformal metric, but the conformal metric will not be analytic or even smooth on $\mathcal{M}^\ast$. In coordinates, $\scri$ is given by ${U_+}_3 {U_-}_3=A$ which also corresponds to  $r=\infty$, and its spacelike nature produces a behaviour near infinity similar to that of a spacelike singularity. Near $\scri^+$, future-directed causal curves are bound to ``go to infinity'' once they enter region IV. Of course, unlike the spacelike singularity in Schwarzschild, no observer or light ray can reach infinity in a finite amount of an affine parameter of these null and timelike geodesics, so no geodesic incompleteness is caused by the dynamics of region IV. 

The geodesic incompleteness comes from the singularity at $r=0$: Radial null geodesics hit the singularity in a finite amount of their affine parameter, so, $\mathcal{M}^\ast$ is geodesically incomplete. However, null and timelike curves can avoid hitting the singularity and go from region II$'$ to region II passing through the ``wormhole''. This indicates that the singularity is timelike. Despite  geodesic incompleteness, the spacetime is timelike geodesically complete as the singularity is repulsive, due to the Reissner-Nordstr{\o}m nature of our spacetime.  To see why, consider for simplicity radial timelike geodesics: A radial curve $\gamma(\tau)=(t(\tau),r(\tau),\omega_0)$, for some constant angular coordinates, is geodesic if
\begin{eqnarray}
\ddot{t}&=&-\dot{t}\dot{r}\frac{f'}{f} \; ;\label{tdoubledot}\\
\ddot{r}&=&-\frac{f'}{2}\left(f\dot{t}^2+\frac{\dot{r}^2}{f} \right) \label{rdoubledot} \; ,
\end{eqnarray}
assuming $f\ne 0$ which is the case near $r=0$, and where dot denote differentiating with respect to $\tau$. In addition, we have $g(\dot{\gamma},\dot{\gamma})=$constant$=E>0$ i.e.
\begin{equation}\label{constantspeedgeodesic}
f \dot{t}^2-\frac{\dot{r}^2}{f}=E \; ,
\end{equation}
So, (\ref{rdoubledot}) becomes
\begin{equation*}
-2\ddot{r}=f'E \; .
\end{equation*}
If we multiple both sides of this equality by $\dot{r}$ then integrate in $\tau$ we obtain
\begin{equation*}
\dot{r}^2+fE=C \; ,
\end{equation*}
where $C$ is the integration constant. We see from (\ref{constantspeedgeodesic}) that $C=f^2 \dot{t}^2$ and hence $C\ge 0$. Thus,
\begin{equation}\label{constrainsonr}
\dot{r}^2=C-fE \; ,
\end{equation}
but $f>0$ on $0<r<r_1$ and in fact $f\rightarrow +\infty$ as $r\rightarrow 0$, which puts a constraint on $r$ preventing it from reaching zero. This means that there must be a turning point in the curve $\gamma$ after which $r$ starts to increase again. So, even objects in free-fall directly (i.e. radially) towards the singularity get ejected to the other region II. Therefore, no timelike geodesic can hit the singularity. The timelike nature of the singularity also means that there are points in the spacetime whose both future and past null cones meet the singularity inside the same region I. Another consequence of this nature is the absence of a (global) Cauchy hypersurface, as there are inextendible timelike curves of arbitrarily small length which start and end at the singularity. For instance, the spacelike hypersurface $\mathscr S$ in figure \ref{fig:OpensConnected1} is a Cauchy hypersurface for regions covered by the domains $B_{k,l-1}$ and $C_{k,l-1}$ for all $k\in \mathbb{Z}$. Yet, there are future-directed and past-directed inextendible timelike curves of $\mathcal{M}^*$ which do not intersect $\mathscr S$. Such curves hit the singularity inside region I and never cross the horizons at $r=r_1$ towards $\mathscr S$, therefore, data in regions I do not depend on data at $\mathscr S$. The hypersurfaces $-\mathscr{H}^-_1 \cup \mathscr{H}^-_1$ and $-\mathscr{H}^+_1 \cup \mathscr{H}^+_1$ bounding regions II and II$'$ in $B_{k,l-1}$ (for all $k$) are said to be Cauchy horizons for the spacelike section $\mathscr S$ (see \cite{hawking_large_1973}). 

%\begin{figure}
%	\centering
%	\includegraphics[scale=1]{Spacelikesections}
%	\caption{\textit{Spacelike Sections}}
%	\label{fig:Spacelikesections}
%\end{figure}

Each point of the diagram shown in figure \ref{fig:OpensConnected1} is a 2-sphere of $\mathcal{M}^*$, or of $\scri^\pm$ which are conformal spacelike 3-hypersurfaces. The segments labelled by $r=0$ where spherical coordinates are singular, are 1-dimensional lines of singular points (of the metric) representing the center of the black hole at different times. We note that the singularity does not touchy null infinity in reality. The corners of the removed squares $S_{k,l}$s, labelled by $i^+$ and $i^-$ and called future and past timelike infinities respectively, are distinct from the segments of $r=0$ because there are plenty of inextendible timelike curves that do not hit the singularities. For example, the timelike curves of constant $r$ in regions III ($r_2<r<r_3$) of the form $\gamma(\tau)=(\tau,C,\omega_0)$ for $\tau\in \R$ in the $(t,r,\omega)$-coordinates never approach the singularity, and one of them is a geodesic, namely when $f'(C)=0$. These future-directed timelike curves originate at $i^-$ and finish at $i^+$. We also note that there are extendible timelike curves that have no end points in the closure of $\mathcal{M}^*$ such as those given locally by $t=0$ in regions II and II$'$.

%Consider a future-directed inextendible timelike geodesic $ \gamma(\tau)$ where $\tau$ is the proper time of the geodesic, and let $$P^\pm=\lim_{\tau\rightarrow \pm \infty}\gamma(\tau) \; .$$ It may happen that both or one of the limit points do not exist, i.e. are not in the closure of $\mathcal{M}^\ast$ in $\R^2 \times \mathcal{S}^2$. Moreover, since no timelike geodesic touches the singularity, $P^\pm$ can't be on the singularity. On the other hand, if $\gamma(\tau)$ belongs to  region IV or IV$'$ for some $\tau$, then either $P^+\in \scri^+$ or $P^-\in \scri^-$. If there is some $A>0$ such that $\gamma(\tau)$ is in some region I for all $\tau>A$, then $P^+$ is in $i^-$. And $P^-$ is in $i^+$ if $\gamma(\tau)$ is in some region I for all $\tau<-A$. Similarly, if the same is true but in region III instead of I, then $P^+$ and $P^-$ exchange roles. This characterizes $i^\pm$ as past and future timelike infinities. It is worth mentioning that there is no analogue to $i_0$ in this spacetime.

\section*{Remark}
We end with a remark about the different number of horizons. The construction carried out in the paper can be easily modified to account for the cases with less number of horizons. The case we treated here is in some sense the most complete. In the case of three horizons, the maximal extension contains all the blocks that can appear in the maximal extensions of the cases with fewer horizons. That is, in the other cases, only some of the blocks I, II, III, IV are present. Also, the conditions on the mass, the charge, and the cosmological constant, for $f$ to have less number of zeros, can be found using arguments similar to those of section \ref{Sec:Zerosoff}.

\section*{Acknowledgement}
The results of this paper and the mentioned decay and conformal scattering results \cite{mokdad_decay_2016,mokdad_conformal_2016}, were obtained during my PhD thesis \cite{mokdad_maxwell_2016}. I would like to thank my thesis advisor Pr. Jean-Philippe Nicolas for his indispensable guidance during the thesis.

\printbibliography[heading=bibintoc]

\end{document}